\documentclass{article}

\usepackage[preprint]{neurips_2024}
\usepackage{tikz}
\usepackage{svg-extract}
\usetikzlibrary{positioning,angles,quotes}
\usepackage{subcaption}
\usepackage{graphicx}
\usepackage{float}
\usepackage{environ}
\usepackage{float}
\usepackage{svg}
\usepackage{amsmath}
\usepackage{comment}
\usepackage{amssymb}
\usepackage{natbib} 
\everypar{\looseness=-1}
\setcitestyle{square,numbers,comma}
\usepackage{adjustbox}
\usepackage{booktabs}
\usepackage{bm}
\usepackage{amsthm}

\newtheorem{theorem}{Theorem}[section]
\newtheorem{lemma}[theorem]{Lemma}
\newtheorem{definition}[theorem]{Definition}
\newtheorem{assumption}[theorem]{Assumption}
\newtheorem{proposition}[theorem]{Proposition}

\newtheorem{prob}[theorem]{Problem}

\usepackage{algorithm,algpseudocode}
\usepackage{hyperref}

\newtheorem{example}[theorem]{Example}

\usepackage{mathdots}
\newcommand{\R}{\mathbb{R}}
\newcommand{\E}{\mathbb{E}}
\newcommand{\B}{\mathcal{B}}
\newcommand{\PP}{\mathbb{P}}
\newcommand{\A}{\mathcal{A}}
\newcommand{\X}{\mathcal{X}}
\newcommand{\St}{\mathcal{S}}
\newcommand{\G}{\mathcal{G}}
\newcommand{\K}{\mathcal{K}}
\newcommand{\ts}{\tilde{s}}
\newcommand{\tpi}{\tilde{\pi}}

\newcommand{\tS}{\tilde{\mathcal{S}}}
\newcommand{\rew}{V_{{r},0}^\theta(\tilde{s}_0)}
\newcommand{\cost}{V_{{c},0}^\theta(\tilde{s}_0)}

\newcommand{\ecost}{\hat{V}_{{c},0}^\theta(\tilde{s}_0)}

\usepackage{environ}
\usepackage{multicol}

\usepackage[utf8]{inputenc} 
\usepackage[T1]{fontenc}    
\usepackage{hyperref}       
\usepackage{cleveref}
\usepackage{url}            
\usepackage{booktabs}       
\usepackage{amsfonts}       
\usepackage{nicefrac}       
\usepackage{microtype}      
\usepackage{xcolor}         
\allowdisplaybreaks[4]

\title{A learning-based approach to stochastic optimal control under reach-avoid constraint}

\author{%
  Tingting Ni \\
  SYCAMORE, EPFL\\
  \texttt{tingting.ni@epfl.ch} \\
   \And
   Maryam Kamgarpour \\
  SYCAMORE, EPFL \\
  maryam.kamgarpour@epfl.ch\\
}

\begin{document}

\maketitle

\begin{abstract}
We develop a model-free approach to optimally control stochastic, Markovian systems subject to a reach-avoid constraint. Specifically, the state trajectory must remain within a safe set while reaching a target set within a finite time horizon. Due to the time-dependent nature of these constraints, we show that, in general, the optimal policy for this constrained stochastic control problem is non-Markovian, which increases the computational complexity. To address this challenge, we apply the state-augmentation technique from \cite{schmid2024joint}, reformulating the problem as a constrained Markov decision process (CMDP) on an extended state space. This transformation allows us to search for a Markovian policy, avoiding the complexity of non-Markovian policies. To learn the optimal policy without a system model, and using only trajectory data, we develop a log-barrier policy gradient approach. We prove that under suitable assumptions, the policy parameters converge to the optimal parameters, while ensuring that the system trajectories satisfy the stochastic reach-avoid constraint with high probability.
\end{abstract}

\section{Introduction}
Many real-world control applications require systems to operate reliably and safely, as failures in fields such as air traffic control~\cite{prandini2008application} or autonomous vehicles \cite{hsu2021safety} can have severe consequences. As complex, safety-critical problems involving stochastic systems emerge, there is a growing need not only to optimize performance but also to enforce constraints that ensure safety, preventing events like collisions, and reachability, which guarantees the successful completion of goals such as cargo delivery. 

To capture these requirements, the reach-avoid specification in a hybrid state space has been considered in \cite{abate2008probabilistic, summers2009probabilistic, summers2011stochastic}, which ensures that the system's states remain within a safe set while also reaching a designated target set. For stochastic systems, this specification can only be satisfied probabilistically, leading to the formulation of a chance constraint on the state trajectory. On the other hand, we aim to optimize performance independently of this specification—for example, by reaching the goal set as quickly as possible. The combination of performance optimization with these constraints results in a chance-constrained stochastic optimal control problem.\looseness-1

Despite the broad applicability of such problems, solving them is challenging due to the time-dependent nature of the constraints, which breaks the problem’s Markovian structure. Consequently, one has to consider history-dependent policies where control actions depend on past trajectories. To find such policies, \cite{thorpe2022data} reformulated the original problem using a dataset of observed trajectories and relaxed it to a linear programming problem. However, the resulting formulation is computationally challenging, even for one-dimensional examples. The authors of~\cite{hahn2019interval} formulated the problem as multiple probabilistic temporal logic objectives on an interval Markov decision process (MDP) and developed a value-iteration algorithm to approximate the Pareto frontier. However, this method is limited to finite state and action spaces and does not address the suboptimality of the resulting policies. Alternatively, the authors of~\cite{wisniewski2023probabilistic} consider using stopping times to apply the evolution equation to devise a linear program for computing the optimal policy. However, this method is also restricted to finite state and action spaces. \looseness-1

Several works~\cite{ono2015chance,paulson2020stochastic,wang2020non} circumvent history-dependence by applying Boole’s inequality to decouple time dependencies, allowing for Markovian policies that rely only on the current state. However, this approach provides a conservative approximation, and optimal tuning of the constraint violation probability at each timestep is also challenging.\looseness-1

Recent works by \cite{schmid2023computing, schmid2024joint} address the non-Markovian structure through state augmentation, reformulating the problem as bilevel optimization within the MDP framework, where the augmented state space is hybrid, as it is a product of continuous and discrete spaces. For policy optimization, \cite{schmid2023computing} uses dynamic programming, while \cite{schmid2024joint} uses linear programming to obtain Markovian policies for each sub-problem. However, these approaches remain computationally challenging in continuous or large discrete state and action spaces, as they rely on the inner steps of dynamic programming or linear programming. Additionally, they require knowledge of the system model.

Building on the state-augmentation techniques from \cite{schmid2023computing,schmid2024joint}, we reformulate the chance-constrained stochastic optimal control problem as a finite-horizon constrained MDP (CMDP). We then focus on a model-free learning approach to approximate the optimal policy based on system trajectory data. Specifically, we consider a setting where we can query the system with a chosen policy to obtain the corresponding system trajectory. Our goal is then to design a learning algorithm that ensures the high-probability satisfaction of the constraints throughout the learning process, while converging to the optimal policy.



To ensure constraint satisfaction during learning for CMDPs, which we refer to as safe exploration, model-free approaches have been developed in the past and demonstrated empirical success in handling complex environments \cite{achiam2017constrained, tessler2018reward, liu2020ipo}. Building on advances in stochastic constrained optimization \cite{usmanova2024log}, \cite{ni2024safe} recently extended the log-barrier policy gradient approach to infinite-horizon discounted CMDPs, showing that it guarantees safe exploration while converging to the optimal policy with high probability. However, the reach-avoid objectives in our problem cannot be framed as a discounted constraint. Unlike infinite-horizon discounted CMDPs, where it suffices to search for the optimal policy within the class of stationary Markov policies, optimal policies in finite-horizon CMDPs are non-stationary. Therefore, extending the theoretical guarantees of the log-barrier policy gradient approach is necessary to address these challenges in the context of finite-horizon CMDPs.

In this paper, we adapt the log-barrier policy gradient method to address the chance-constrained stochastic optimal control problem. Our contributions are as follows:
\begin{enumerate}
    \item While past work has alluded to the necessity of optimizing over history-dependent policies by considering all possible policies \cite{schmid2023computing,schmid2024joint}, we demonstrate concretely in Example \ref{example} that generally, the optimal policy for reach-avoid constrained stochastic optimal control cannot be Markovian.
    \item Inspired by \cite{schmid2023computing, schmid2024joint}, we show that the chance-constrained stochastic optimal control problem can be reformulated as a CMDP on an augmented state space (see Theorem \ref{thm_pro_pro}).
    \item We derive a learning-based approach to approximate the solution with a provable guarantee while ensuring a high probability of constraint satisfaction during the learning process (see Theorem \ref{main}).
\end{enumerate}
The rest of the paper is organized as follows. In Section \ref{section:problemformulation}, we formulate the chance-constrained stochastic optimal control problem and demonstrate that the general optimal policy for such a problem is non-Markovian. In Section \ref{sec_r_CMDP}, we prove that this problem can be equivalently reformulated as a CMDP. In Section \ref{sec_LGPG}, we introduce Algorithm \ref{alg:cap} which employs the log-barrier policy gradient method. In Section \ref{sec_conan}, we provide proof sketches and intuition behind the convergence and safe exploration guarantees of Algorithm \ref{alg:cap}, with full technical details provided in the appendix to maintain clarity. We conclude with a summary and outlook in Section \ref{sec_c_fd}.

\subsection{Notation}
Let $\mathbb{N}$ and $\R$ denote the set of natural and real numbers, respectively, and for any $m \in \mathbb{N}$, define $[m]:= \{0, 1, \dots, m\}$. For any topological space $\X$, we denote by $\B(\X)$ the Borel $\sigma$-algebra on that space. For any subset $\mathcal{Y}\subset \B(\X)$, the indicator function of the set $\mathcal{Y}$ is denoted by  $\mathbf{1}_{\mathcal{Y}}(\cdot)$. Here, $\|\cdot\|$ denotes the Euclidean $\ell_2$-norm for vectors and the operator norm for matrices, respectively. For any two matrices $A,B\in\R^{n\times n }$, $A \succeq B$ indicates that the matrix {$A-B$} is positive semi-definite. We denote the Moore–Penrose inverse of $A$ by $A^\dagger$. A function $f:\X\to \R^n$, is said to be {$M$}-smooth if {$f(x)\le f(y)+\langle\nabla f(y),x-y\rangle+\frac{M}{2}\left\|x-y\right\|^2$} holds {$ \forall x,y\in \X$}, and is {$L$}-Lipschitz continuous if {$\left|f(x)-f(y)\right|\le L\left\|x-y\right\|$} holds { $\forall x,y\in \X$}.
We denote the logical conjunction and disjunction by $\wedge$ and $\vee$, respectively.
\section{Problem formulation}\label{section:problemformulation}
We consider a finite-horizon Markov decision process (MDP) defined by the tuple {$ \{\mathcal{S},\mathcal{A},H, P,r\}$}. The state space {$\mathcal{S}$} and the action space {$\mathcal{A}$} are Borel subsets of complete and separable metric spaces, each equipped with the $\sigma$-algebras $\B(\St)$ and $\B(\A)$, respectively. The time horizon is denoted by $H\in\mathbb{N}$. The transition kernel {${P}(s'|s,a)$} represents the probability density of transitioning from state $s$ to state $s'$ under action $a$. The reward $r:=\{\left\{r_t\right\}_{t\in[H-1]},r_H\}$ is a set of $H+1$ measurable functions, where $r_t:\mathcal{S}\times \mathcal{A} \rightarrow [0,1]$ represents the stage reward for all $t\in[H-1]$, and {$r_H:\mathcal{S} \rightarrow [0,1]$} is the terminal reward.

We further denote by $\mathbb{H}_t$ the space of all histories up to timestep $t$ with $\mathbb{H}_0 = \St$ at the initial timestep and $ \mathbb{H}_t = (\St \times \A)^{t} \times \St$ at the following timesteps $t\in\{1,\dots,H\}$. An element of $h_t\in\mathbb{H}_t$ is called a trajectory and is of the form $h_t:=\{(s_0,a_0),\dots,(s_{t-1},a_{t-1}),s_t\}$. A policy $\pi$ is a set of Borel-measurable stochastic kernels $\{\pi_t\}_{t\in[H-1]}$, where each $\pi_t(\cdot|h_t)$ assigns a probability measure on the set $\B(\A)$ for each $h_t\in\mathbb{H}_t$. A Markov policy $\pi$ is a policy, where each $\pi_t$ depends only on the current state, namely $\pi_t(\cdot|h_t)=\pi_t{(\cdot|s_t)}$ for all $h_t\in\mathbb{H}_t$ and $t\in[H-1]$. We denote the sets of all policies and Markov policies by $\Pi$ and $\Pi_{\text{M}}$, respectively. Note that $\Pi_M \subset \Pi$.


Let $(\mathbb{H}_{H}, \B(\mathbb{H}_{H}))$ be a measurable space. By \cite[Proposition V.1.1]{neveu1965mathematical}, there exists a unique probability measure $P^{\pi}$ on $\B(\mathbb{H}_{H})$ for a given initial state  $s_0\in \St$, policy $\pi\in\Pi$ and the transition kernel $P$. We denote by $\E^{\pi}$ the expectation with respect to the probability measure  $P^{\pi}$. Then, the expected cumulative reward at timestep $t\in[H-1]$ for a given state $s$ is defined as follows:
\begin{align}
V_{r, t}^\pi(s) := \E^{\pi} \left[ \sum_{i=t}^{H-1} r_{i}(s_i, a_i) + r_H(s_H) \mid s_t =s\right].
\end{align}
Since \( P^{\pi} \) is unique for any given policy $\pi$, it follows that the expected cumulative reward \( V_{r, t}^\pi(s) \) is uniquely determined as well.

Next, we consider a safe set and a goal set. Specifically, let \(\K \in \B(\mathcal{S})\) denote the safe set and \(\G \in \B(\mathcal{S})\) the goal set, with \(\G \subset \K\). A given trajectory $\tau\in\mathbb{H}_{H}$ satisfies the reach-avoid property if the following condition is met:
\[
\mathcal{C}_{\K,\G} := \{\tau \in \mathbb{H}_{H} \mid \exists t \in [H], s_t \in \G, s_i \in \K\setminus\G, \forall i \in [t-1] \}.
\]
Note that the reach-avoid property encompasses the invariance and reachability properties, as shown in \cite{tkachev2013quantitative}.

In this paper, we aim to find a policy $\pi^*\in\Pi$ that maximizes the expected cumulative reward $V_{r,0}^\pi(s_0)$, while simultaneously ensuring that the probability of satisfying reach-avoid property exceeds a specified threshold $\delta$. We formalize this so-called chance-constrained stochastic optimal control as follows:
\begin{prob}\label{Original_problem}
Find {$\pi^*\in\Pi$} such that it solves the problem below:
\begin{align}
\sup_{\pi \in \Pi} V_{r,0}^\pi(s_0)\quad \text{s.t.} \quad P_{s_0}^{\pi}(\mathcal{C}_{\K,\G}) \ge \delta, 
\end{align}
where $\delta \in [0,1)$ is the risk tolerance parameter. 
\end{prob}
Our objective is to solve the problem without knowledge of the transition dynamics, relying instead on access to system trajectories. Specifically, for any given policy, we can observe the corresponding trajectory. Thus, our goal is twofold: on the one hand, we aim to ensure that the queried policies satisfy the constraints during learning (safe exploration), and on the other hand, we seek to ensure that the queried policies converge to the optimal policy.

To ensure that Problem \ref{Original_problem} is well-posed, we make the following assumption.
\begin{assumption}\label{Assumption_continuity}
\begin{enumerate}
\item The action space $\mathcal{A}$ is compact.
\item For any $s \in \mathcal{S}$ and $s' \in \mathcal{B}(\mathcal{S})$, the transition kernel  $P(s' | s, a)$  and the stage reward  $r_t(s, a)$, for all $t\in [H-1]$, are continuous with respect to $a$.
\item There exists a policy $\pi_0$ and a positive constant $\nu_s$ such that  $P_{s_0}^{\pi_0}(\mathcal{C}_{\K,\G}) \ge \delta+\nu_s$.
\end{enumerate}
\end{assumption}
Point (1) of Assumption \ref{Assumption_continuity} ensures the compactness of the action space \( \mathcal{A} \), while Point (2) requires continuity of both the transition kernel and stage reward functions. These conditions are automatically satisfied if both the action space $\mathcal{A}$ and state space $\mathcal{S}$ are finite, assuming the discrete topology on $\mathcal{A}$ and $\mathcal{S}$, respectively. Lastly, Point (3), known as Slater’s condition, ensures strict feasibility of Problem \ref{Original_problem}. This condition can be verified by computing the maximum probability of satisfying the reach-avoid property through linear programming over state-action occupancy measures. Alternatively, this probability can be approximated using dynamic programming with a discretization procedure on the continuous state space \cite{tkachev2013quantitative}.
As will be demonstrated in Theorem \ref{thm_pro_pro}, Assumption \ref{Assumption_continuity} ensures the existence of an optimal policy attaining the supremum in Problem~\ref{Original_problem}.
\subsection{Necessity of non-Markovian policies}\label{sec_n_o_p}
In this section, we illustrate the need to optimize over the set of non-Markovian policies $\Pi$ for Problem \ref{Original_problem}. Since the chance constraint cannot be equivalently decomposed into individual stage-wise constraints, the optimal policy at timestep $t$ may depend on whether the history $h_t$ violates or satisfies the reach-avoid property. As illustrated in the following example, if the reach-avoid property is already violated at timestep \( t \), the agent may prioritize reward maximization. Conversely, if the reach-avoid property has not been violated, the agent may focus on ensuring its satisfaction.
\begin{example}\label{example}
    Consider the MDP shown in Figure~\ref{graph_mdp} which has six states $\St=\{s^i\}_{i\in[5]}$, two actions \(\mathcal{A} = \{a^i\}_{i\in[1]}\), and time horizon $H=3$. The transition dynamics are defined as follows: with equal probability $s^0$ transitions either to $s^1$ or $s^2$; and both $s^1$ and $s^2$ transition deterministically to $s^3$. The state $s^3$  transitions deterministically to $s^4$ and $s^5$ under actions $a^0$ and $a^1$, respectively. The stage rewards $\{r_t\}_{t\in[2]}$ are zero, and the terminal reward $r_3(s) = \mathbf{1}_{s=s^5}$, for all $s\in\St$. We define the safe set as $\K = \{s^0, s^1, s^3, s^4\}$ and the goal set as $\G = \{s^4\}$. In Figure \ref{graph_mdp}, we mark $\G$ in red, and $\K \setminus \G$ in green. The problem is to solve the following chance-constrained stochastic optimal control problem with initial state $s^0$:
    \begin{align}
\sup_{\pi \in \Pi} V_{r,0}^\pi(s^0)\quad \text{s.t.} \quad P_{s^0}^{\pi}(\mathcal{C}_{\K,\G}) \ge 0.4.
\end{align}
The strictly feasible policy is to always take action \( a_0 \), ensuring that the reach-avoid property is satisfied with probability 0.5. Therefore, the problem satisfies Assumption~\ref{Assumption_continuity}. 

We note that at state \( s^3 \), taking action \( a^0 \) satisfies the reach property but does not contribute to the reward, whereas taking action \( a^1 \) violates the reach property while optimizing for the reward. To balance reward optimization and constraint satisfaction, the optimal policy at state \( s^3 \) and time \( t = 2 \) is given by \( \pi_2(a^0 \mid s^0, s^1, s^3) = 0.8 \) and \( \pi_2(a^1 \mid s^0, s^2, s^3) = 1 \), resulting in a reward of \( V_{r,0}^\pi(s^0) = 0.6 \). However, if restricted to Markov policies, the optimal policy for $s^3$ at timestep $t=2$ is $\pi_2^M(a^0|s^3)=0.8$, with reward $V_{r,0}^{\pi^M}(s^0)=0.1$. Clearly, the reward $V_{r,0}^\pi(s^0)=0.6$ is larger than reward $V_{r,0}^{\pi^M}(s^0)=0.1$.
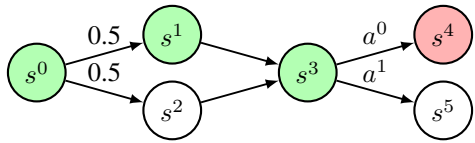
\begin{figure}[h]
\centering
\begin{center}
\begin{tikzpicture}[auto,node distance=1cm,>=latex,font=]
\tikzstyle{round}=[thick,draw=black,circle]
    \node[round, fill=green!30] (s0) {$s^0$};
    \node[round, fill=green!30, right=of s0, yshift=0.5cm] (s1) {$s^1$};
    \node[round, right=of s0, yshift=-0.5cm] (s2) {$s^2$};
    \node[round, fill=green!30, right=of s1, yshift=-0.5cm] (s3) {$s^3$};
    \node[round, fill=red!30, right=of s3, yshift=0.5cm] (s4) {$s^4$};
    \node[round, right=of s3, yshift=-0.5cm] (s5) {$s^5$};
    \path[->, draw=black,solid,line width=0.25mm,fill=black] (s0) -- node[above] {0.5} (s1);
    \path[->, draw=black,solid,line width=0.25mm,fill=black] (s0) -- node[above] {0.5} (s2);
    \path[->, draw=black,solid,line width=0.25mm,fill=black] (s1) -- node[above] {} (s3);
    \path[->, draw=black,solid,line width=0.25mm,fill=black] (s2) -- node[above] {} (s3);
    \path[->, draw=black,solid,line width=0.25mm,fill=black] (s3) -- node[above] {$a^0$} (s4);
    \path[->, draw=black,solid,line width=0.25mm,fill=black] (s3) -- node[above] {$a^1$} (s5);
\end{tikzpicture}
\end{center}
\caption{Six-states MDP with safe set $\K = \{s^0, s^1, s^3, s^4\}$ and goal set $\G = \{s^4\}$.}
\label{graph_mdp}
\end{figure}
\vspace{-0.5cm}
\end{example}
In the above example, any optimal policy is non-Markovian. This observation also holds for invariance and reachability constraints when considering only the safe set $\K$. 

Given that optimizing over the set $\Pi$ is challenging, recent work by \cite{schmid2024joint, tkachev2013quantitative} addresses this issue through a state-augmentation technique that transforms the non-Markovian constraint into one with a Markovian structure. In the next section, we apply this technique to convert the original problem into a constrained Markov decision process (CMDP). Under Assumption \ref{Assumption_continuity}, this CMDP has an optimal Markov policy that corresponds to the optimal non-Markovian policy in the original MDP.
\section{Formulation as a constrained Markov decision process}\label{sec_r_CMDP}
To address the non-Markovian structure of the reach-avoid constraint, \cite{schmid2024joint, tkachev2013quantitative, wang2022solving} introduced a state-augmentation technique that reformulates the reach-avoid probability as an objective defined over the final augmented state, allowing for an equivalent and tractable Markovian structure. Here, we adopt the approach from ~\cite{schmid2024joint}, as it allows us to represent the constraint within the CMDP framework, whose well-established theory offers a computationally efficient solution.

Specifically, we consider a CMDP defined by the tuple \(\{\tilde{\mathcal{S}}, \mathcal{A},H, \allowbreak\tilde{P}, {r}, {c}\}\), where $\tS:=\St\times \mathcal{Y}$ is the augmented state space with $\mathcal{Y}:=\{0,1,2\}$. Here, $0$, $1$, and $2$ capture the events of reaching the unsafe set $\St\setminus\K$, the safe but not goal set $\K\setminus\G$, and the goal set $\G$, respectively. Each element of $\tS$ is of the form $\ts=(s,y)$, where $s\in\St$ and $y\in\mathcal{Y}$. The time horizon $H$ and action space $\A$ remain the same as in the original MDP.

Let us define the set of valid transitions among the auxiliary states as $\mathcal{Y}_\text{valid}$, specified by the conditions:
\begin{small}
\begin{align}
(y_t\! = y_{t+1}\!  =\!  0) \! \vee \! (y_t \! = \! y_{t+1} \! = 2) \! \vee\!  (y_t \!= \!1 \!\wedge\! y_{t+1} \!= \!\mathbf{1}_{s_{t+1} \in \mathcal{K}\setminus\G}\! + \!2\cdot \mathbf{1}_{s_{t+1} \in \mathcal{G}}).
\end{align}
\end{small}

\noindent This condition implies that states $0$ and $2$ are absorbing within the state space $\mathcal{Y}$, as reaching the unsafe set or goal set indicates that the reach-avoid objective has been violated or satisfied, respectively. State $1$, on the other hand, can transition to states $0$, $1$, or $2$, representing possible transitions to the unsafe set, a safe but not goal set, or the goal set, respectively. 

Then, the transition kernel $\tilde{P}$ is defined as:
\begin{align}
    &\tilde{P}((s_{t+1},y_{t+1})|(s_t,y_t),a_t) \\
    =& \begin{cases}
        {P}({s}_{t+1}|{s}_{t},a_t), &\text{if }y_{t+1}=f(s_{t+1},y_t) \wedge (y_t,y_{t+1})\in \mathcal{Y}_\text{valid}\\
         0, &\text{otherwise,}
    \end{cases}
\end{align}
where $f:\St\times\mathcal{Y}\to\mathcal{Y}$  is given by:
\begin{align}
    f(s,y) = \begin{cases}
        2, \text{ if } y = 2 \vee (y = 1 \wedge s\in\G),\\
        1, \text{ if } y = 1 \wedge s\in\K\setminus\G,\\
        0, \text{ otherwise. } 
\end{cases}\label{eq_dyna}
\end{align}
 To summarize, the transition process unfolds in two steps: Firstly, the state $s_t$ transitions to $s_{t+1}$ based on $P(\cdot|s_t, a_t)$ as in the MDP setting. Secondly, if the auxiliary state $y_{t+1}$ equals $f(s_{t+1},y_t)$ and $(y_t,y_{t+1})$ lies in the set of valid transitions $\mathcal{Y}_\text{valid}$, the transition dynamic remains the same as $P(s_{t+1}|s_t, a_t)$; otherwise, it is set to zero. 

For this new MDP, the reward \( r \) is inherited from the original MDP. The stage-wise constraints are set to zero, and the terminal constraint is defined as \( c(\tilde{s}) = \mathbf{1}_{y=2} \), for all \( \tilde{s} \in \tilde{\mathcal{S}} \). Policies are defined as before. Let \( \tilde{\Pi} \) and \( \tilde{\Pi}_{\text{M}} \) denote the sets of all policies and Markov policies on the augmented state space \( \tilde{\mathcal{S}} \), respectively. Note that \( \tilde{\Pi}_{\text{M}} \subset \tilde{\Pi} \), and any policy in \( \Pi \) can be transformed into a policy in \( \tilde{\Pi} \) by setting \( \tilde{\pi}(\cdot \mid s,y) := \pi(\cdot \mid s) \) for any \( y \in \mathcal{Y} \) and \( \pi \in \Pi \).
Given the initial state $\ts_0:=(s_0,\mathbf{1}_{s_0 \in \mathcal{K}\setminus\G} + 2\cdot \mathbf{1}_{s_0 \in \mathcal{G}})$, a policy $\tilde{\pi}\in\tilde{\Pi}_{\text{M}}$, and the transition kernel $\tilde{P}$, a trajectory of the form $\{(s_0,y_0),a_0,\dots,(s_H,y_H)\}$ is generated. The design of the transition kernel $\tilde{P}$ ensures that the reach-avoid property ~\cite[Lemma 3.1]{schmid2024joint} can be verified via the auxiliary states $\{y_t\}_{t\in[H]}$. Namely, reach-avoid is satisfied along the trajectory if and only if $y_H = 2$.

Let ${P}^{\tilde{\pi}}$ denote the probability measure over the trajectory space $(\tilde{\mathcal{S}} \times \mathcal{A})^H \times \tilde{\mathcal{S}}$, which is uniquely defined by \cite[Proposition V.1.1]{neveu1965mathematical}. We denote the corresponding expectation by $\mathbb{E}^{\tilde{\pi}}$. Then, the expected cumulative reward at timestep $t\in[H-1]$ for a given state $\tilde{s}$ is given by:
\begin{align}
    V_{{r}, t}^{\tilde{\pi}}(\tilde{s}) := \E^{\tilde{\pi}} \left[ \sum_{i=t}^{H-1} r_{i}(s_i, a_i) + r_H(s_H) \mid \ts_t =\ts\right],
\end{align}
and the expected cumulative constraint value is given by:
\begin{align}
    V_{{c}, t}^{\tilde{\pi}}(\tilde{s}) :=\E^{\tilde{\pi}}\left[{c}(\tilde{s}_H) \mid \tilde{s}_t =\tilde{s}\right].
\end{align}
Unlike the chance constraint \( P_{s_0}^{\pi}(\mathcal{C}_{\K, \G}) \), the reformulated problem imposes only a final-stage constraint, thereby eliminating any need for coupling constraints across the time horizon. Given the expected reward \( V_{r,0}^{\tilde{\pi}}( \tilde{s}_0) \) and expected constraint value \( V_{c,0}^{\tilde{\pi}}(\tilde{s}_0) \), we formulate the CMDP problem as follows:

\begin{prob}\label{AMDP_problem}
Find \( \tilde{\pi} \in \tilde{\Pi} \) such that it solves the problem below:
\[
\sup_{\tilde{\pi} \in \tilde{\Pi}} V_{r,0}^{\tilde{\pi}}( \tilde{s}_0) \quad \text{s.t.} \quad V_{c,0}^{\tilde{\pi}}( \tilde{s}_0) \ge \delta,
\]
where \( \delta \in [0, 1) \) is the risk tolerance parameter.
\end{prob}
Next, we demonstrate that solving Problem \ref{Original_problem} is equivalent to solving Problem \ref{AMDP_problem}.
\begin{theorem}\label{thm_pro_pro}
    Under Assumption \ref{Assumption_continuity}, there exists an optimal Markovian policy {\small\( \tilde{\pi}_M^* \)} for Problem \ref{AMDP_problem} that corresponds to an optimal non-Markovian policy for Problem \ref{Original_problem}. Specifically, we construct a mapping \( g: \tilde{\Pi}_M \to \Pi \) such that
\begin{align}
    V_{r,t}^{\pi^*}(s_0) = V_{r,t}^{g(\tilde{\pi}_M^*)}(s_0) \quad \text{s.t.} \quad P_{s_0}^{g(\tilde{\pi}_M^*)} \ge \delta,
\end{align}
where the mapping \( g \) is defined in the proof below (see equation \eqref{eq_defg}) and $\pi^*$ is an optimal policy for Problem \ref{Original_problem}.
\end{theorem}
\begin{proof}
Under Assumption \ref{Assumption_continuity}, there exists a strictly feasible policy \( \pi_0 \) for Problem \ref{Original_problem}, which is also strictly feasible for Problem~\ref{AMDP_problem} since \( \Pi \subset \tilde{\Pi} \). Additionally, due to the design of the transition kernel \( \tilde{P} \), \( \tilde{P}(\tilde{s}' | \tilde{s}, a) \) is continuous with respect to \( a \) for any \( \tilde{s} \in \tilde{\mathcal{S}} \) and \( \tilde{s}' \in \mathcal{B}(\tilde{\mathcal{S}}) \). Since Problem \ref{AMDP_problem} can be formulated as a linear constrained optimization over the state-action occupancy measure rather than over policy parameters \( \pi \),\footnote{Here, we define the state-action occupancy measure \( \{d_t^{\tpi}(\tilde{s},a)\}_{t \in [H-1]}\), for any \((\tilde{s}, a) \in \tilde{\St} \times \A \), as $d_t^{\tpi}(\tilde{s}, a)= \mathbb{E}^{\tpi}\left[\mathbf{1}_{(\tilde{s}_t,a_t):=(\tilde{s},a)}| \ts_0\right]$.} the existence of a strictly feasible policy ensures strong duality for Problem \ref{AMDP_problem} \cite[Corollary 2]{piunovskii1994control}. Moreover, the compactness of the action space and the continuity of both the transition kernel and stage reward functions guarantee the attainment of the supremum for Problem \ref{AMDP_problem}. Consequently, an optimal Markovian policy \( \tilde{\pi}_M^* \) exists for Problem \ref{AMDP_problem} under Assumption~\ref{Assumption_continuity}~\cite[Theorem 3]{piunovskii1994control}.

To show that an optimal Markovian policy \( \tilde{\pi}_M^* \) for Problem \ref{AMDP_problem} results in an optimal non-Markovian policy for Problem \ref{Original_problem}, we start by mapping a history $h_t=\{(s_0,a_0),\dots,(s_{t-1},a_{t-1}),s_t\}$ to an auxiliary state \( b_t \) using the transition dynamics from \eqref{eq_dyna}. Specifically, we calculate \( b_t \) recursively as follows:
\begin{align}
    b_0 \!:=\! \mathbf{1}_{s_0 \in \mathcal{K}\setminus\G} + 2 \cdot \mathbf{1}_{s_0 \in \mathcal{G}}, \text{and }b_t \!:=\! f(s_{t}, b_{t-1}), \forall t\in\{1,\dots,H-1\}.
\end{align}
This allows us to associate a Markovian policy \( \pi_t(\cdot|s_t,b_t) \), defined over the augmented space, with a non-Markovian policy \( \pi_t(\cdot|h_t) \) based on the full history. Consequently, we construct a projection map \( g : \tilde{\Pi}_{\text{M}} \to \Pi \) as follows: for any policy \( \pi \in \Pi \),
\begin{align}
   g^{-1}\left( \{\pi_t(a_t|h_t)\}_{t\in[H-1]}\right)=\{\tpi_t(a_t|s_t,b_t)\}_{t\in[H-1]}.\label{eq_defg}
\end{align}
Since \( g(\tilde{\pi}) \) produces the same distribution over states as \( \pi \), we have the following relations for any \( t \in [H-1] \), \( s \in \mathcal{S} \), and \( \tilde{\pi} \in \tilde{\Pi}_{\text{M}} \):
\begin{equation}
    \begin{aligned}
    &V_{{r}, t}^{\tilde{\pi}}((s, \mathbf{1}_{s \in \mathcal{K}\setminus\G} + 2 \cdot \mathbf{1}_{s \in \mathcal{G}})) = V_{{r}, t}^{g(\tilde{\pi})}(s) , \,
    V_{{c}, 0}^{\tilde{\pi}}(\tilde{s}_0) = P_{s_0}^{g(\tilde{\pi})}(\mathcal{C}_{\mathcal{K},\mathcal{G}}). \label{eq_equation}
    \end{aligned}
\end{equation}
Therefore,
\begin{align}
 &V_{{r}, t}^{\tilde{\pi}_M^*}(\tilde{s}_0) = V_{{r}, t}^{g(\tilde{\pi}_M^*)}(s_0) \le V_{{r}, t}^{\pi^*}(s_0), \,V_{{c}, 0}^{\tilde{\pi}_M^*}(\tilde{s}_0) = P_{s_0}^{g(\tilde{\pi}_M^*)}(\mathcal{C}_{\mathcal{K},\mathcal{G}}) \ge \delta.
\end{align}
Since the optimal policy \( \pi^* \) for Problem \ref{Original_problem} is also an element of \( \tilde{\Pi} \), we have
\begin{align}
&V_{{r}, t}^{\pi^*}(s_0) = V_{{r}, t}^{\pi^*}(\tilde{s}_0) \le V_{{r}, t}^{\tilde{\pi}_M^*}(\tilde{s}_0).
\end{align}
Consequently, \( g(\tilde{\pi}_M^*) \) is an optimal policy for Problem \ref{Original_problem}.
\end{proof}
Leveraging Theorem \ref{thm_pro_pro}, we reduce our objective to finding the optimal Markovian policy for Problem \ref{AMDP_problem}, while ensuring \emph{safe exploration}, as defined below.
\begin{definition}\label{definition_safe_exploration}
An algorithm generating a sequence $\{\tilde{\pi}(i)\}_{i \in [I]}$  satisfies safe exploration if \( \,V_{c,0}^{\tilde{\pi}(i)}(\tilde{s}_0) \ge \delta \) holds for all \( i \in [I] \).
\end{definition}
Observe that, as shown in \eqref{eq_equation}, if \( \tilde{\pi} \) is feasible for Problem~\ref{AMDP_problem}, then \( g(\tilde{\pi}) \) is feasible for Problem \ref{Original_problem}. Thus, if the algorithm ensures safe exploration for Problem \ref{AMDP_problem}, it also guarantees safe exploration for the original problem.

To efficiently search for the optimal policy in large or continuous CMDPs, we employ policy gradient methods with a finite-dimensional parameterization of the policy \cite{agarwal2021theory}. This approach enables us to reformulate the search for Markovian policies in Problem~\ref{AMDP_problem} as a constrained optimization problem over a finite-dimensional parameter space, as shown in Problem \ref{Paramterized_problem}. After this reformulation, we present a learning algorithm that updates the policy parameters such that each resulting policy satisfies the constraint, namely \(V_{c,0}^{\tilde{\pi}(i)}(\tilde{s}_0) \ge \delta\) with high probability. Additionally, under suitable assumptions, the parameters converge to the optimal policy.

\section{Log barrier policy gradient approach} \label{sec_LGPG}
We begin by parameterizing the Markovian policies as \(\{\tpi_t^{\theta_t}\}_{t \in [H-1]}\), where $\theta_t$ denotes the parameter of the policy at timestep $t$. With this finite-dimensional parameterization, we can cast Problem \ref{AMDP_problem} as a finite-dimensional constrained optimization problem.

For example, in the case of a continuous action space within the interval \([a, b]\), a truncated Gaussian distribution \cite{mihoc2003fisher}  can be employed and is defined as:  
\begin{align}
    \pi_\theta(a|s) = {\phi( \frac{a - \theta(s)}{\sigma} )}/{\sigma \bigl( \Phi( \frac{b - \theta(s)}{\sigma} ) - \Phi( \frac{a - \theta(s)}{\sigma} ) \bigr)},
\end{align}  
where \( \theta(s) \) represents the mean, \( \sigma^2 \) is the variance, and \( \phi \) and \( \Phi \) denote the probability density function and cumulative distribution function of the standard normal distribution, respectively. For finite action spaces, a log-linear policy is commonly used \cite{yuan2023linear}, defined as:  
\begin{align}
    \pi_\theta(a|s) = {\exp(f_\theta(s,a))}/{\sum_{a' \in \mathcal{A}} \exp(f_\theta(s,a'))},
\end{align}  
where \( f_\theta(s,a) \) is the feature function. Above, if  \( f_\theta(s,a) = \theta(s,a) \), the parametrization is referred to as the softmax, and if \( f_\theta(s,a) \) is parameterized using a neural network, it is referred to as neural softmax. For hybrid systems as is the case in our Problem \ref{AMDP_problem}, the policy parameterization can be constructed as a product of two separate parameterizations \cite{delalleau2019discrete}, each defined for the continuous and discrete components of the state space. For additional choices of policy parameterization methods, see \cite{rajeswaran2017towards}.

For simplicity, let \(\theta := \{\theta_t\}_{t\in[H-1]}\) and \(\tpi^\theta := \{\tpi_t^{\theta_t}\}_{t \in [H-1]}\). We abbreviate \(\E^{\tpi^{\theta}}\), \(V_{{r},t}^{\tpi^{\theta}}(\tilde{s})\), and \(V_{{c},t}^{\tpi^{\theta}}(\tilde{s})\) as \(\E^{\theta}\), \(V_{{r},t}^{\theta}(\tilde{s})\), and \(V_{{c},t}^{\theta}(\tilde{s})\), respectively. With this parameterization, we reformulate Problem~\ref{AMDP_problem} as follows:
\begin{prob}\label{Paramterized_problem}
    Find \(\theta \in \R^d\) such that it solves the problem below:
    \begin{align}
        \sup_{\theta \in \R^d} V_{{r},0}^{\theta}(\tilde{s}_0) \quad \text{s.t.} \quad V_{{c},0}^{\theta}(\tilde{s}_0) \ge \delta,
    \end{align}
    where \(\delta \in [0,1)\). The feasible set is denoted by {\small\(\Theta := \{\theta \mid V_{{c},0}^{\theta}(\tilde{s}_0) \ge \delta\}\)}.
\end{prob}
Observe that the above is a finite-dimensional constrained optimization problem. However, several challenges arise: (i) the objective function, constraints, and their gradients are unknown; (ii) the problem is non-convex; and (iii) interacting with the system directly requires satisfying the constraints. To address challenge (iii), we algorithmically apply the log-barrier method. This interior-point approach ensures the feasibility of the iterates $\theta(i)$. To address challenges (ii) and (i), we apply gradient ascent on the log-barrier, where we estimate the gradients using a Monte Carlo approach by sampling trajectories from the system. Finally, to ensure convergence to a stationary point of the log-barrier function, we show that our problem enjoys a so-called gradient dominance property under suitable assumptions.

Inspired by interior-point algorithms, we ensure safe exploration by maximizing an unconstrained log-barrier surrogate of Problem~\ref{Paramterized_problem}.
\begin{align}
    B_\eta^{\theta} := V_{{r},0}^{\theta}(\tilde{s}_0) + \eta \log(V_{{c},0}^{\theta}(\tilde{s}_0) - \delta), \label{equation_log_barrier}
\end{align}
where \(\eta > 0\). We maximize this objective using stochastic gradient ascent with updates given by
\begin{align}
    \label{equation_updates}
    \theta(i+1) = \theta(i) + \gamma(i) \hat{\nabla}_\theta B_\eta^{\theta(i)},
\end{align}
where \(\gamma(i)\) is the stepsize and \(\hat{\nabla}_\theta B_\eta^{\theta(i)}\) estimates the true gradient \(\nabla_\theta B_\eta^{\theta(i)}\), calculated as:
\begin{align}
    \nabla_\theta B_\eta^{\theta} = \nabla_\theta V_{{r},0}^{\theta}(\ts_0) + \eta {\nabla_\theta V_{{c},0}^{\theta}(\ts_0)}/({V_{{c},0}^{\theta}(\ts_0) - \delta}). \label{equation_loggradient}
\end{align}
The log-barrier method begins with a safe policy \(\pi^{\theta(0)}\) and updates the policy parameters iteratively by \eqref{equation_updates}. With a carefully chosen stepsize \(\gamma(i)\), the log-barrier algorithm ensures safe exploration while converging to the optimal policy parameters (see Section~\ref{sec_conan}). In the following two sections, we will explain how to estimate the gradient using system trajectories, as well as how to choose the stepsize to ensure safe exploration and convergence.

\subsection{Estimating the log barrier gradient}\label{section_estimatinggradient}
To estimate the gradient \(\nabla_\theta B_\eta^{\theta}\) without direct access to the transition dynamics \({P}\), we estimate it using sample trajectories from the MDP as follows.

Given a policy parameter $\theta$,  we apply $\pi^\theta$ to the system to obtain \(n\) trajectories, denoted by {\small\(\{\tau_i\}_{i \in [n-1]}\)}. For each trajectory \(\tau_i\), we collect the reward and constraint values, with the sampled data denoted as 
\[\{\tilde{s}_t^i, a_t^i, {r}_{t}({s}_t^i, a_t^i)\}_{t\in [H-1]} \cup \{{r}_H({s}^i_H), {c}(\tilde{s}^i_H)\}.\]
The estimator for \(\cost\), denoted by \(\ecost\), is computed as the sample mean over these trajectories:
\[\hat V_{c,0}^\theta(\ts_0)  :=\frac{1}{n}\sum_{i\in[n-1]}{c}(\tilde{s}^i_H).\]
It is known that the estimators of \( \nabla_\theta \rew \) and \( \nabla_\theta \cost \), denoted by \( \hat{\nabla}_\theta \rew \) and \( \hat{\nabla}_\theta \ecost \), respectively, can be calculated as follows \cite{klein2024beyond}:
\begin{small}
\begin{align}
&\hat\nabla_\theta \rew\!:= \! \frac{1}{n}\!\!\sum_{i\in[n-1]}\!\sum_{t\in[H-1]}\!\!\!\nabla_\theta\log \tpi^{\theta_t}_t\!(a^i_t|\ts^i_t) (\!\sum_{m=t}^{H-1} \!{r}_m(\!{s}^i_m, a^i_m\!) \!+\! {r}_H({s}^i_H\!) \!),\\
&\hat\nabla_\theta \cost:=\frac{1}{n}\sum_{i\in[n-1]}\sum_{t\in[H-1]}\nabla_\theta\log \tpi^{\theta_t}_t(a^i_t|\ts^i_t)){c}(\tilde{s}^i_H).
\end{align}
\end{small}

To bound the error of these estimations, we rely on the assumption about Lipschitz and smoothness properties of the policy.
\begin{assumption}
\label{assumption_smooth_policy} 
There exist non-negative constants $M_g$ and $M_h$ such that for all $t\in[H-1]$, $\theta\in\Theta$, and $(\ts,a)\in\tS\times \A$, the gradient and Hessian of $\,\log \tpi^{\theta_t}_t$ satisfy 
\begin{align}
    \|\nabla_{\theta_t}\log \tpi^{\theta_t}_t(a|\ts)\|\le M_g,\,\|\nabla_{\theta_t}^2\log \tpi^{\theta_t}_t(a|\ts)\|\le M_h.
\end{align}
\end{assumption}
This assumption is satisfied by log-linear policies with bounded feature vectors \cite[Section 6.1.1]{agarwal2021theory} and by truncated Gaussian parameterizations with bounded variance \cite{mihoc2003fisher}. Under Assumption~\ref{assumption_smooth_policy}, we establish the smoothness and Lipschitz continuity of both {\small{\(\rew\)}} and {\small{\(\cost\)}}, enabling sub-Gaussian tail bounds on their value and gradient estimators, as shown below.

\begin{proposition}\label{pro_smooth}
Let Assumption \ref{assumption_smooth_policy} hold. Then the following properties are satisfied for all {$\theta\in\Theta$}:

1. For any $k\in\{r,c\}$, the function {\small$V_{k,0}^\theta(\ts_0)$} is $L_k$-Lipschitz continuous and $M_k$-smooth, where 
$L_r:=M_gH^{\frac{3}{2}}$, $L_{c}:=H^{\frac{1}{2}}M_g$, $M_r:=M_hH+M_g^2H^2$ and $M_{c}:=M_h+M_g^2H$

2. The estimators {\small{\(\ecost\), \(\hat\nabla_\theta \cost\)}}, and {\small{\(\hat\nabla_\theta \rew\)}} are unbiased.

3. For any $\beta\in(0,1)$, we have
\begin{align}
    &\PP\left(\left|\ecost- \cost\right|\le \sigma_{c}^0(n)\sqrt{\ln{\frac{2}{\beta}}}\right)\ge 1-\beta,
\end{align}
where {$\sigma_{c}^0(n):={1}/{\sqrt{2n}}$}.

4. For any $\beta\in(0,1)$, the following holds for any $k\in\{r,c\}$ by setting $n\ge8\ln{\frac{\exp{\frac{1}{4}}}{\beta}}$:
\begin{align}
    &\PP\left(\left\|\hat\nabla_\theta V_{k,0}^\theta(\ts_0) -\nabla_\theta V_{k,0}^\theta(\ts_0)\right\|\le  \sigma_{k}^1(n)\sqrt{\ln{\frac{\exp{\frac{1}{4}}}{\beta}}}\right)\ge 1-\beta,
\end{align}
where {$\sigma_r^1(n):={2\sqrt{2}L_r}/{\sqrt{n}}$} and {$\sigma_{c}^1(n):={2\sqrt{2}L_c}/{\sqrt{n}}$}.
\end{proposition}
The above results extend \cite[Proposition 3.2]{ni2024safe} from infinite-horizon discounted CMDPs to finite-horizon CMDPs. For completeness, we provide the proof in \cite[Appendix A.1]{ni2024learning}. Using these estimators, we approximate \( \nabla_\theta B_\eta^{\theta} \) by
\begin{align}
    \hat{\nabla}_\theta B_\eta^{\theta} := \hat{\nabla}_\theta \rew + \eta {\hat{\nabla}_\theta \cost}/(\ecost - \delta). \label{eq_estimatelog}
\end{align}
Particularly, the concentration bounds ensure that {\small\( \hat{\nabla}_\theta B_\eta^{\theta} \)} remains close to the true gradient {\small\( \nabla_\theta B_\eta^{\theta} \)}. This property allows us to derive high-probability guarantees for both safe exploration and algorithmic convergence.

\subsection{Learning algorithm}
\label{sec_LA}
To implement the gradient ascent in \eqref{equation_updates}, the stepsize should be small enough to ensure that the resulting policy \( \tilde{\pi}^{\theta(i)} \) satisfies the constraint, namely \( V_{c,0}^{{\theta(i)}}(\tilde{s}_0) \ge \delta \) with high probability. At the same time, the stepsize should also be large enough to allow the iterates to make progress. Both of these considerations can be addressed given a smoothness constant, which bounds the growth of the gradient of the function. However, our log barrier function \( B_\eta^{\theta} \) is not globally smooth, as its gradient diverges to infinity near the boundary of the feasible domain. To address this, \cite{hinder2023worst} examines the local smoothness of \( B_\eta^{\theta} \), where its gradient exhibits bounded growth within a local region. However, this constant depends on the gradients of the objective and constraints, which are unknown in our model-free setting.

Inspired by the stochastic optimization framework in \cite{usmanova2024log}, we estimate the smoothness constant based on our stochastic function and gradient estimates, as outlined in Proposition \ref{pro_smooth}. In particular, we construct the lower confidence bound of \( V_{c,0}^{\theta(i)}(\ts_0) - \delta \) and the upper confidence bound of {\small\(|\langle \nabla_\theta V_{c,0}^{\theta(i)}(\ts_0), {\nabla_\theta B_\eta^{\theta(i)}}/{\|\nabla_\theta B_\eta^{\theta(i)}\|} \rangle| \)} using Proposition \ref{pro_smooth}, denoted by \( \underline{\alpha}(i) \) and \( \overline{\beta}(i) \), respectively. We calculate \( \underline{\alpha}(i) \) and \( \overline{\beta}(i) \) as
\begin{align}
    \underline{\alpha}(i) &:= \hat V_{c,0}^{\theta(i)}(\ts_0) - \delta - \sigma_c^0(n)\sqrt{\ln\frac{2}{\beta}}, \\
    \overline{\beta}(i) &:= \left|\left\langle \hat \nabla_\theta V_{c,0}^{\theta(i)}(\ts_0), \frac{\hat\nabla_\theta B_\eta^{\theta(i)}}{\|\hat\nabla_\theta B_\eta^{\theta(i)}\|}\right\rangle\right| + \sigma_c^1(n)\sqrt{\ln\frac{\exp{\frac{1}{4}}}{\beta}}.
\end{align}
This enables us to estimate the local smoothness constant \( \hat{M}(i) \), calculated as
\begin{align}
    \hat{M}(i) := M_r + {10 M_c \eta}/{\underline{\alpha}(i)} + 8 \eta {(\overline{\beta}(i))^2}/{\left(\underline{\alpha}(i)\right)^2}. \label{eq_estM}
\end{align}
By the concentration bounds for {\small\( V_{c,0}^{\theta(i)}(\ts_0) \)} and its gradient {\small\( \nabla_\theta V_{c,0}^{\theta(i)}\)}, \( \hat{M}(i) \) remains close to the true local smoothness constant \( M(i) \).

To prevent overshooting and ensure iterations stay within the valid region in which $\hat{M}(i)$ is a valid bound on the gradient growth of the log barrier, we choose the stepsize \( \gamma(i) \) as
\begin{align}
    \gamma(i) = \min\left\{ {1}/{\hat{M}(i)}, {\underline{\alpha}(i)}/\left({( \sqrt{M_c \underline{\alpha}(i)} + 2 |\overline{\beta}(i)| ) \| \hat{\nabla}_\theta B_\eta^{\theta(i)} \|}\right) \right\}. \label{gammat}
\end{align}
The second term ensures that the stepsize remains within the valid region around the current iterate \( \theta(i) \).

Combining the above steps, Algorithm \ref{alg:cap} can be summarized as: if the norm of the estimated gradient exceeds \({\eta}/{2}\), the algorithm performs stochastic gradient ascent using the gradient estimators defined in~\eqref{eq_estimatelog} and the stepsize defined in~\eqref{gammat}. Otherwise, the algorithm terminates.
\begin{algorithm}
\caption{}\label{alg:cap}
\begin{algorithmic}[1]
\State \textbf{Input}: Batch size $n$, iteration number $I$, $\beta\in(0,1)$, $\eta>0$.
\For{$i=0,1,\dots,I-1$}
\State Compute {$\hat\nabla_\theta B_\eta^{\theta(i)}$} using Eq.\eqref{eq_estimatelog}.
\If{{$\|\hat\nabla_\theta B_\eta^{\theta(i)}\| \le \frac{\eta}{2}$}}
\State Break and return $\theta({\text{break}}):=\theta(i)$.
\EndIf
\State {$\theta(i+1)=\theta(i)+\gamma(i) \hat\nabla_\theta B_\eta^{\theta(i)}$}, where $\gamma(i)$ is defined in Eq.\eqref{gammat}.
\EndFor
\State Return $\theta({\text{out}})$, which can be either $\theta({\text{break}})$ or $\theta(I)$.
\end{algorithmic}
\end{algorithm}
\subsection{Convergence analysis}\label{sec_conan}
In this section, we show that the Algorithm \ref{alg:cap} ensures safe exploration and finds an \( \tilde{\mathcal{O}}(\varepsilon) \)-optimal policy with high probability, using \( \tilde{\mathcal{O}}(\varepsilon^{-6}) \) trajectories sampled from the system, as detailed in Theorem \ref{main}. \footnote{The notation $\Tilde{\mathcal{O}}(\cdot)$ hides the $\log(\frac{1}{\varepsilon})$ term.}

To establish safe exploration guarantees, we adopt the extended Mangasarian-Fromovitz constraint qualification (MFCQ) assumption, which ensures that for any point sufficiently close to the boundary (within a distance of $p$), there exists a direction which is \( \nabla_\theta \cost \) pointing away from the boundary and back toward the interior of the feasible region (by $\ell$).
\begin{assumption}[Extended MFCQ]\label{emf}
    There exist positive constants \( p, \ell \) such that for all \( \theta \in \{\theta \in \Theta \mid 0 < \cost - \delta \le p \} \), we have \( \|\nabla_\theta \cost \| \ge \ell \).
\end{assumption}
This assumption is adopted from the log barrier for constrained stochastic optimization setting \cite{usmanova2024log}.

Given the model-free setting, to ensure sufficient exploration, past work on policy gradient methods \cite{yuan2022general, agarwal2021theory, fazel2018global, liu2020improved} has relied on the Fisher non-degeneracy assumption. Furthermore, \cite{ni2024safe} relaxed this assumption to the relaxed Fisher non-degeneracy condition (see Assumption \ref{fn}). We define the state-action occupancy measure at timestep \( t \) as \( d_t^{\theta}(\ts,a) := \mathbb{E}^{\theta}[\mathbf{1}_{(\ts_t,a_t)=(\ts,a)}|\ts_0 ] \), and the Fisher information matrix as
\[
F_t^\theta := \mathbb{E}_{(\ts,a) \sim d_t^\theta} [\nabla_{\theta_t} \log \tpi_t^\theta(a|\ts)( \nabla_{\theta_t} \log \tpi_t^\theta(a|\ts) )^T ].
\]
Under Assumption 4.2, the Fisher information matrix is guaranteed to be well-defined and positive semi-definite \cite{kakade2001natural}.
\begin{assumption}[Relaxed Fisher non-degeneracy]\label{fn}
    For all \( \theta \in \Theta \) and \( t \in [H-1] \), there exists a positive constant \( \mu_F\) such that the smallest non-zero eigenvalue of \( F_t^\theta \) is lower bounded by \( \mu_F \).
\end{assumption}

Assumption \ref{fn} is satisfied by properties of the transition kernel that ensure sufficient exploration of each state-action pair, combined with a suitable policy parameterization. For instance, softmax, log-linear, and neural softmax parameterizations satisfy Assumption \ref{fn} by constraining $\Theta$ to a compact set \cite{ni2024safe}. Additionally, truncated Gaussian parameterization also meets this assumption~\cite{mihoc2003fisher}.

Next, we introduce the concept of transfer error, which quantifies how well the policy set \( \{\tpi_\theta\}_{\theta \in \mathbb{R}^d} \) approximates the set of all stochastic Markovian policies. We begin by defining the state-action reward and constraint functions as
\begin{align}
&Q_{r,t}^{\theta}(\tilde{s}, a) := \mathbb{E}^{\theta} [ \sum_{i=t}^{H-1} r_i(s_i, a_i) + r_H(s_H) \mid (\ts_t, a_t) = (\ts, a) ],\\
&Q_{c,t}^{\theta}(\tilde{s}, a) := \mathbb{E}^{\theta} \left[ c(\tilde{s}_H) \mid (\ts_t, a_t) = (\ts, a) \right],
\end{align}
respectively. For any function \( k \in \{r, c\} \), we define the advantage function as \( A_{k,t}^{\theta}(\tilde{s}, a) = Q_{k,t}^{\theta}(\tilde{s}, a) - V_{k,t}^{\theta}(\tilde{s}) \) and the transfer error as
\[
L(\mu_{k,t}, \theta, d_t^{\tpi^*_M}) := \mathbb{E}_{(s,a) \sim d_t^{\tpi^*_M}} [ (A_{k,t}^\theta(\ts,a) - \mu_{k,t}^T \nabla_{\theta_t} \log \tpi_t^\theta(a|\ts) )^2 ],
\]
where \( \mu_{k,t} = (F_t^\theta)^\dagger \mathbb{E}_{(\ts,a) \sim d_t^\theta} [ \nabla_{\theta_t} \log \tpi_t^\theta(a|\ts)^T A_{k,t}^\theta(\ts, a) ] \).

We define the transfer error of a policy class as
\begin{align}
    \varepsilon_{bias}:=\sup\left\{ L(\mu_{k,t}, \theta, d_t^{\tpi^*_M}) \mid \forall\theta \in \Theta,\,t \in [H-1], \, k \in \{r,c\}\right\}.
\end{align}
To ensure suboptimality of the output policy returned by Algorithm~\ref{alg:cap}, as shown in Theorem \ref{main}, it is essential to minimize \( \varepsilon_{bias} \), which is achieved by parameterizations such as softmax with \( \varepsilon_{bias} = 0 \), and rich log-linear, or neural softmax with small transfer error~\cite{ni2024safe, wang2019neural}.

Having established all the preceding assumptions, we now present the main theorem.
\begin{theorem}\label{main}
    Let Assumptions~\ref{Assumption_continuity}, ~\ref{assumption_smooth_policy}, ~\ref{emf}, and ~\ref{fn} hold. For any {\small\(\beta\in(0,1)\)}, by setting {\small$\eta\le\min\{p,\nu_s,({8H^{\frac{1}{2}}M_g}/{(C\mu_F)})^{\frac{1}{2}}\}$} and {\small\(n = \mathcal{O}(\eta^{-4}\ln\frac{1}{\beta})\)}, the following holds after \(I\) iterations of Algorithm \ref{alg:cap}:
    \begin{enumerate}
        \item \(\PP\left(\min_{i\in[I]} V_{0,c}^{\theta(i)}(\ts_0) \ge c\eta\right) \ge 1 - I\beta\).
        \item With a probability of at least \(1 - I\beta\), the output policy \(\tpi^{\theta({\text{out}})}\) satisfies
            \begin{align}
                V_{r,0}^{\tpi^*}(\ts_0) - V_{r,0}^{\theta(\text{out})}(\ts_0) 
                \le& \left(1 - \frac{C\mu_F\eta^2}{8H^{\frac{1}{2}}M_g}\right)^{I} \left(V_{r,0}^{\tpi^*}(\ts_0) - V_{r,0}^{\theta(0)}(\ts_0)\right) \\
                & + \mathcal{O}({\sqrt{\varepsilon_{bias}}}) + \tilde{\mathcal{O}}(\eta),
            \end{align}
            where $c:=\frac{\ell}{24L_c}$ and \(C := \frac{c}{2{L_r}^2(1 + \frac{1}{c})\max\left\{4 + \frac{5{M_r}c}{{L_r}^2}, 1 + \sqrt{\frac{{M_r}c\eta}{4{L_r}^2}}\right\}}\)
    \end{enumerate}
\end{theorem}
The proof of Theorem \ref{main} consists of two parts: safe exploration analysis and convergence analysis. In the safe exploration analysis, we establish a lower bound on the distance of the iterates from the boundary in \cite[Lemma A.1]{ni2024learning}. In the convergence analysis, we demonstrate the gradient dominance property of the log barrier function in \cite[Lemma A.2]{ni2024learning}, which enables us to bound the gap between \( V_{r,0}^{\tpi^*}(\ts_0) \) and \( \rew \) by the norm of the gradient \( \nabla_\theta B_\eta^\theta \). Furthermore, the gradient ascent method ensures convergence to a stationary point of \( B_\eta^\theta \). Leveraging these results, we complete the proof of Theorem \ref{main}, with a detailed proof provided in \cite[Appendix A.2]{ni2024learning}.

By setting \( \eta = \varepsilon \) and \( I = \tilde{\mathcal{O}}(\varepsilon^{-2}) \) in the above theorem, Algorithm \ref{alg:cap} ensures, with high probability, that the original constraint, namely 
\( P_{s_0}^{g(\tpi^{\theta(i)})}(\mathcal{C}_{\K,\G}) \ge \delta \), is satisfied throughout the learning process and converges to the optimal set of parameters within the policy class using \( \tilde{\mathcal{O}}(\varepsilon^{-6}) \) trajectories sampled from the system. Additionally, the returned policy \( g(\pi^{\theta(\text{out})} )\) has a bounded transfer error of \( \mathcal{O}(\sqrt{\varepsilon_{{bias}}}) \) relative to the optimal policy for the original problem (see Problem \ref{Original_problem}).

This result extends the findings of \cite{ni2024safe} from infinite-horizon discounted CMDPs to finite-horizon CMDPs, demonstrating that the sample complexity remains consistent across both settings and advancing our understanding of algorithmic efficiency.

\paragraph{Discussion: } 
The chance-constrained stochastic optimal control model considered in this work optimizes rewards for all trajectories, even those that violate the reach-avoid constraints. This can lead to the controller deliberately failing to meet the constraints in order to achieve a higher reward, as shown in Example \ref{example}. One possible way to avoid this situation is to only collect rewards if the reach-avoid property is satisfied. For instance, we can define the value function as:
\[
V_{r, t}^\pi(s) = \mathbb{E}^{\pi} \left[\left( \sum_{i=t}^{H-1} r_{i}(s_i, a_i) + r_H(s_H)\right)\mathbf{1}_{(s_0, a_0, \dots, s_H)\in\mathcal{C}_{\K, \G}} \mid s_t = s\right].
\]
However, optimizing such a reward subject to the reach-avoid constraint raises the following questions: Is it sufficient to optimize over Markovian policies for global optimality? If not, does there exist an augmented space where optimizing over Markovian policies in this augmented space is sufficient? We will leave these as future work.
\section{Conclusions and future directions} \label{sec_c_fd}
We considered a state-augmentation technique to convert the chance-constrained stochastic optimal control problem into a finite-horizon CMDP. This transformation ensures that the optimal Markovian policy derived from the CMDP aligns with the optimal non-Markovian policy of the original problem. To find this policy while ensuring safe exploration, we employed the log-barrier algorithm proposed in \cite{usmanova2024log}. Through rigorous analysis, we proved that this algorithm guarantees safe exploration and convergence to the optimal policy.

A promising direction for future work is the design of policy parameterizations with meaningful feature functions. Another important avenue is to conduct extensive numerical simulations to better understand the potential and limitations of the algorithm, both in benchmark and real-world settings.
\section{Acknowledgments}
This research is gratefully supported by the Swiss National Science Foundation (SNSF).
\bibliographystyle{ACM-Reference-Format}
\bibliography{ref}


\begin{thebibliography}{36}


\ifx \showCODEN    \undefined \def \showCODEN     #1{\unskip}     \fi
\ifx \showDOI      \undefined \def \showDOI       #1{#1}\fi
\ifx \showISBNx    \undefined \def \showISBNx     #1{\unskip}     \fi
\ifx \showISBNxiii \undefined \def \showISBNxiii  #1{\unskip}     \fi
\ifx \showISSN     \undefined \def \showISSN      #1{\unskip}     \fi
\ifx \showLCCN     \undefined \def \showLCCN      #1{\unskip}     \fi
\ifx \shownote     \undefined \def \shownote      #1{#1}          \fi
\ifx \showarticletitle \undefined \def \showarticletitle #1{#1}   \fi
\ifx \showURL      \undefined \def \showURL       {\relax}        \fi
\providecommand\bibfield[2]{#2}
\providecommand\bibinfo[2]{#2}
\providecommand\natexlab[1]{#1}
\providecommand\showeprint[2][]{arXiv:#2}

\bibitem[Abate et~al\mbox{.}(2008)]%
        {abate2008probabilistic}
\bibfield{author}{\bibinfo{person}{Alessandro Abate}, \bibinfo{person}{Maria Prandini}, \bibinfo{person}{John Lygeros}, {and} \bibinfo{person}{Shankar Sastry}.} \bibinfo{year}{2008}\natexlab{}.
\newblock \showarticletitle{Probabilistic reachability and safety for controlled discrete time stochastic hybrid systems}.
\newblock \bibinfo{journal}{\emph{Automatica}} \bibinfo{volume}{44}, \bibinfo{number}{11} (\bibinfo{year}{2008}), \bibinfo{pages}{2724--2734}.
\newblock


\bibitem[Achiam et~al\mbox{.}(2017)]%
        {achiam2017constrained}
\bibfield{author}{\bibinfo{person}{Joshua Achiam}, \bibinfo{person}{David Held}, \bibinfo{person}{Aviv Tamar}, {and} \bibinfo{person}{Pieter Abbeel}.} \bibinfo{year}{2017}\natexlab{}.
\newblock \showarticletitle{Constrained policy optimization}. In \bibinfo{booktitle}{\emph{International conference on machine learning}}. PMLR, \bibinfo{pages}{22--31}.
\newblock


\bibitem[Agarwal et~al\mbox{.}(2021)]%
        {agarwal2021theory}
\bibfield{author}{\bibinfo{person}{Alekh Agarwal}, \bibinfo{person}{Sham~M Kakade}, \bibinfo{person}{Jason~D Lee}, {and} \bibinfo{person}{Gaurav Mahajan}.} \bibinfo{year}{2021}\natexlab{}.
\newblock \showarticletitle{On the theory of policy gradient methods: Optimality, approximation, and distribution shift}.
\newblock \bibinfo{journal}{\emph{Journal of Machine Learning Research}} \bibinfo{volume}{22}, \bibinfo{number}{98} (\bibinfo{year}{2021}), \bibinfo{pages}{1--76}.
\newblock


\bibitem[Delalleau et~al\mbox{.}(2019)]%
        {delalleau2019discrete}
\bibfield{author}{\bibinfo{person}{Olivier Delalleau}, \bibinfo{person}{Maxim Peter}, \bibinfo{person}{Eloi Alonso}, {and} \bibinfo{person}{Adrien Logut}.} \bibinfo{year}{2019}\natexlab{}.
\newblock \showarticletitle{Discrete and continuous action representation for practical rl in video games}.
\newblock \bibinfo{journal}{\emph{arXiv preprint arXiv:1912.11077}} (\bibinfo{year}{2019}).
\newblock


\bibitem[Fazel et~al\mbox{.}(2018)]%
        {fazel2018global}
\bibfield{author}{\bibinfo{person}{Maryam Fazel}, \bibinfo{person}{Rong Ge}, \bibinfo{person}{Sham Kakade}, {and} \bibinfo{person}{Mehran Mesbahi}.} \bibinfo{year}{2018}\natexlab{}.
\newblock \showarticletitle{Global convergence of policy gradient methods for the linear quadratic regulator}. In \bibinfo{booktitle}{\emph{International conference on machine learning}}. PMLR, \bibinfo{pages}{1467--1476}.
\newblock


\bibitem[Hahn et~al\mbox{.}(2019)]%
        {hahn2019interval}
\bibfield{author}{\bibinfo{person}{Ernst~Moritz Hahn}, \bibinfo{person}{Vahid Hashemi}, \bibinfo{person}{Holger Hermanns}, \bibinfo{person}{Morteza Lahijanian}, {and} \bibinfo{person}{Andrea Turrini}.} \bibinfo{year}{2019}\natexlab{}.
\newblock \showarticletitle{Interval Markov decision processes with multiple objectives: from robust strategies to Pareto curves}.
\newblock \bibinfo{journal}{\emph{ACM Transactions on Modeling and Computer Simulation (TOMACS)}} \bibinfo{volume}{29}, \bibinfo{number}{4} (\bibinfo{year}{2019}), \bibinfo{pages}{1--31}.
\newblock


\bibitem[Hinder and Ye(2023)]%
        {hinder2023worst}
\bibfield{author}{\bibinfo{person}{Oliver Hinder} {and} \bibinfo{person}{Yinyu Ye}.} \bibinfo{year}{2023}\natexlab{}.
\newblock \showarticletitle{Worst-case iteration bounds for log barrier methods on problems with nonconvex constraints}.
\newblock \bibinfo{journal}{\emph{Mathematics of Operations Research}} (\bibinfo{year}{2023}).
\newblock


\bibitem[Hsu et~al\mbox{.}(2021)]%
        {hsu2021safety}
\bibfield{author}{\bibinfo{person}{Kai-Chieh Hsu}, \bibinfo{person}{Vicen{\c{c}} Rubies-Royo}, \bibinfo{person}{Claire~J Tomlin}, {and} \bibinfo{person}{Jaime~F Fisac}.} \bibinfo{year}{2021}\natexlab{}.
\newblock \showarticletitle{Safety and liveness guarantees through reach-avoid reinforcement learning}.
\newblock \bibinfo{journal}{\emph{arXiv preprint arXiv:2112.12288}} (\bibinfo{year}{2021}).
\newblock


\bibitem[Kakade(2001)]%
        {kakade2001natural}
\bibfield{author}{\bibinfo{person}{Sham~M Kakade}.} \bibinfo{year}{2001}\natexlab{}.
\newblock \showarticletitle{A natural policy gradient}.
\newblock \bibinfo{journal}{\emph{Advances in neural information processing systems}}  \bibinfo{volume}{14} (\bibinfo{year}{2001}).
\newblock


\bibitem[Klein et~al\mbox{.}(2024)]%
        {klein2024beyond}
\bibfield{author}{\bibinfo{person}{Sara Klein}, \bibinfo{person}{Simon Weissmann}, {and} \bibinfo{person}{Leif D{\"o}ring}.} \bibinfo{year}{2024}\natexlab{}.
\newblock \showarticletitle{Beyond Stationarity: Convergence Analysis of Stochastic Softmax Policy Gradient Methods}. In \bibinfo{booktitle}{\emph{The Twelfth International Conference on Learning Representations}}.
\newblock


\bibitem[Kohler and Lucchi(2017)]%
        {kohler2017sub}
\bibfield{author}{\bibinfo{person}{Jonas~Moritz Kohler} {and} \bibinfo{person}{Aurelien Lucchi}.} \bibinfo{year}{2017}\natexlab{}.
\newblock \showarticletitle{Sub-sampled cubic regularization for non-convex optimization}. In \bibinfo{booktitle}{\emph{International Conference on Machine Learning}}. PMLR, \bibinfo{pages}{1895--1904}.
\newblock


\bibitem[Liu et~al\mbox{.}(2020a)]%
        {liu2020ipo}
\bibfield{author}{\bibinfo{person}{Yongshuai Liu}, \bibinfo{person}{Jiaxin Ding}, {and} \bibinfo{person}{Xin Liu}.} \bibinfo{year}{2020}\natexlab{a}.
\newblock \showarticletitle{Ipo: Interior-point policy optimization under constraints}. In \bibinfo{booktitle}{\emph{Proceedings of the AAAI conference on artificial intelligence}}. \bibinfo{pages}{4940--4947}.
\newblock


\bibitem[Liu et~al\mbox{.}(2020b)]%
        {liu2020improved}
\bibfield{author}{\bibinfo{person}{Yanli Liu}, \bibinfo{person}{Kaiqing Zhang}, \bibinfo{person}{Tamer Basar}, {and} \bibinfo{person}{Wotao Yin}.} \bibinfo{year}{2020}\natexlab{b}.
\newblock \showarticletitle{An improved analysis of (variance-reduced) policy gradient and natural policy gradient methods}.
\newblock \bibinfo{journal}{\emph{Advances in Neural Information Processing Systems}}  \bibinfo{volume}{33} (\bibinfo{year}{2020}), \bibinfo{pages}{7624--7636}.
\newblock


\bibitem[Mihoc and F{\u{a}}tu(2003)]%
        {mihoc2003fisher}
\bibfield{author}{\bibinfo{person}{Ion Mihoc} {and} \bibinfo{person}{Cristina~Ioana F{\u{a}}tu}.} \bibinfo{year}{2003}\natexlab{}.
\newblock \showarticletitle{Fisher's information measures and truncated normal distributions (II)}.
\newblock \bibinfo{journal}{\emph{Revue d'analyse num{\'e}rique et de th{\'e}orie de l'approximation}} \bibinfo{volume}{32}, \bibinfo{number}{2} (\bibinfo{year}{2003}), \bibinfo{pages}{177--186}.
\newblock


\bibitem[Neveu et~al\mbox{.}(1965)]%
        {neveu1965mathematical}
\bibfield{author}{\bibinfo{person}{Jacques Neveu}, \bibinfo{person}{Amiel Feinstein}, {and} \bibinfo{person}{Robert Fortet}.} \bibinfo{year}{1965}\natexlab{}.
\newblock \showarticletitle{Mathematical foundations of the calculus of probability}.
\newblock  (\bibinfo{year}{1965}).
\newblock


\bibitem[Ni and Kamgarpour(2023)]%
        {ni2024safe}
\bibfield{author}{\bibinfo{person}{Tingting Ni} {and} \bibinfo{person}{Maryam Kamgarpour}.} \bibinfo{year}{2023}\natexlab{}.
\newblock \showarticletitle{A safe exploration approach to constrained Markov decision processes}.
\newblock \bibinfo{journal}{\emph{arXiv preprint arXiv:2312.00561}} (\bibinfo{year}{2023}).
\newblock


\bibitem[Ni and Kamgarpour(2024)]%
        {ni2024learning}
\bibfield{author}{\bibinfo{person}{Tingting Ni} {and} \bibinfo{person}{Maryam Kamgarpour}.} \bibinfo{year}{2024}\natexlab{}.
\newblock \showarticletitle{A learning-based approach to stochastic optimal control under reach-avoid constraint}.
\newblock \bibinfo{journal}{\emph{arXiv preprint arXiv:2412.16561}} (\bibinfo{year}{2024}).
\newblock


\bibitem[Ono et~al\mbox{.}(2015)]%
        {ono2015chance}
\bibfield{author}{\bibinfo{person}{Masahiro Ono}, \bibinfo{person}{Marco Pavone}, \bibinfo{person}{Yoshiaki Kuwata}, {and} \bibinfo{person}{J Balaram}.} \bibinfo{year}{2015}\natexlab{}.
\newblock \showarticletitle{Chance-constrained dynamic programming with application to risk-aware robotic space exploration}.
\newblock \bibinfo{journal}{\emph{Autonomous Robots}}  \bibinfo{volume}{39} (\bibinfo{year}{2015}), \bibinfo{pages}{555--571}.
\newblock


\bibitem[Paulson et~al\mbox{.}(2020)]%
        {paulson2020stochastic}
\bibfield{author}{\bibinfo{person}{Joel~A Paulson}, \bibinfo{person}{Edward~A Buehler}, \bibinfo{person}{Richard~D Braatz}, {and} \bibinfo{person}{Ali Mesbah}.} \bibinfo{year}{2020}\natexlab{}.
\newblock \showarticletitle{Stochastic model predictive control with joint chance constraints}.
\newblock \bibinfo{journal}{\emph{Internat. J. Control}} \bibinfo{volume}{93}, \bibinfo{number}{1} (\bibinfo{year}{2020}), \bibinfo{pages}{126--139}.
\newblock


\bibitem[Piunovskii(1994)]%
        {piunovskii1994control}
\bibfield{author}{\bibinfo{person}{AB Piunovskii}.} \bibinfo{year}{1994}\natexlab{}.
\newblock \showarticletitle{Control of random sequences in problems with constraints}.
\newblock \bibinfo{journal}{\emph{Theory of Probability \& Its Applications}} \bibinfo{volume}{38}, \bibinfo{number}{4} (\bibinfo{year}{1994}), \bibinfo{pages}{751--762}.
\newblock


\bibitem[Prandini and Hu(2008)]%
        {prandini2008application}
\bibfield{author}{\bibinfo{person}{Maria Prandini} {and} \bibinfo{person}{Jianghai Hu}.} \bibinfo{year}{2008}\natexlab{}.
\newblock \showarticletitle{Application of reachability analysis for stochastic hybrid systems to aircraft conflict prediction}. In \bibinfo{booktitle}{\emph{2008 47th IEEE conference on decision and control}}. IEEE, \bibinfo{pages}{4036--4041}.
\newblock


\bibitem[Rajeswaran et~al\mbox{.}(2017)]%
        {rajeswaran2017towards}
\bibfield{author}{\bibinfo{person}{Aravind Rajeswaran}, \bibinfo{person}{Kendall Lowrey}, \bibinfo{person}{Emanuel~V Todorov}, {and} \bibinfo{person}{Sham~M Kakade}.} \bibinfo{year}{2017}\natexlab{}.
\newblock \showarticletitle{Towards generalization and simplicity in continuous control}.
\newblock \bibinfo{journal}{\emph{Advances in neural information processing systems}}  \bibinfo{volume}{30} (\bibinfo{year}{2017}).
\newblock


\bibitem[Schmid et~al\mbox{.}(2023)]%
        {schmid2023computing}
\bibfield{author}{\bibinfo{person}{Niklas Schmid}, \bibinfo{person}{Marta Fochesato}, \bibinfo{person}{Sarah~HQ Li}, \bibinfo{person}{Tobias Sutter}, {and} \bibinfo{person}{John Lygeros}.} \bibinfo{year}{2023}\natexlab{}.
\newblock \showarticletitle{Computing optimal joint chance constrained control policies}.
\newblock \bibinfo{journal}{\emph{arXiv preprint arXiv:2312.10495}} (\bibinfo{year}{2023}).
\newblock


\bibitem[Schmid et~al\mbox{.}(2024)]%
        {schmid2024joint}
\bibfield{author}{\bibinfo{person}{Niklas Schmid}, \bibinfo{person}{Marta Fochesato}, \bibinfo{person}{Tobias Sutter}, {and} \bibinfo{person}{John Lygeros}.} \bibinfo{year}{2024}\natexlab{}.
\newblock \showarticletitle{Joint chance constrained optimal control via linear programming}.
\newblock \bibinfo{journal}{\emph{IEEE Control Systems Letters}} (\bibinfo{year}{2024}).
\newblock


\bibitem[Summers et~al\mbox{.}(2011)]%
        {summers2011stochastic}
\bibfield{author}{\bibinfo{person}{Sean Summers}, \bibinfo{person}{Maryam Kamgarpour}, \bibinfo{person}{John Lygeros}, {and} \bibinfo{person}{Claire Tomlin}.} \bibinfo{year}{2011}\natexlab{}.
\newblock \showarticletitle{A stochastic reach-avoid problem with random obstacles}. In \bibinfo{booktitle}{\emph{Proceedings of the 14th international conference on Hybrid systems: computation and control}}. \bibinfo{pages}{251--260}.
\newblock


\bibitem[Summers and Lygeros(2009)]%
        {summers2009probabilistic}
\bibfield{author}{\bibinfo{person}{Sean Summers} {and} \bibinfo{person}{John Lygeros}.} \bibinfo{year}{2009}\natexlab{}.
\newblock \showarticletitle{A probabilistic reach-avoid problem for controlled discrete time stochastic hybrid systems}.
\newblock \bibinfo{journal}{\emph{IFAC Proceedings Volumes}} \bibinfo{volume}{42}, \bibinfo{number}{17} (\bibinfo{year}{2009}), \bibinfo{pages}{150--155}.
\newblock


\bibitem[Tessler et~al\mbox{.}(2019)]%
        {tessler2018reward}
\bibfield{author}{\bibinfo{person}{Chen Tessler}, \bibinfo{person}{Daniel~J. Mankowitz}, {and} \bibinfo{person}{Shie Mannor}.} \bibinfo{year}{2019}\natexlab{}.
\newblock \showarticletitle{Reward constrained policy optimization}. In \bibinfo{booktitle}{\emph{International Conference on Learning Representations}}.
\newblock


\bibitem[Thorpe et~al\mbox{.}(2022)]%
        {thorpe2022data}
\bibfield{author}{\bibinfo{person}{Adam Thorpe}, \bibinfo{person}{Thomas Lew}, \bibinfo{person}{Meeko Oishi}, {and} \bibinfo{person}{Marco Pavone}.} \bibinfo{year}{2022}\natexlab{}.
\newblock \showarticletitle{Data-driven chance constrained control using kernel distribution embeddings}. In \bibinfo{booktitle}{\emph{Learning for Dynamics and Control Conference}}. PMLR, \bibinfo{pages}{790--802}.
\newblock


\bibitem[Tkachev et~al\mbox{.}(2013)]%
        {tkachev2013quantitative}
\bibfield{author}{\bibinfo{person}{Ilya Tkachev}, \bibinfo{person}{Alexandru Mereacre}, \bibinfo{person}{Joost-Pieter Katoen}, {and} \bibinfo{person}{Alessandro Abate}.} \bibinfo{year}{2013}\natexlab{}.
\newblock \showarticletitle{Quantitative automata-based controller synthesis for non-autonomous stochastic hybrid systems}. In \bibinfo{booktitle}{\emph{Proceedings of the 16th international conference on Hybrid systems: computation and control}}. \bibinfo{pages}{293--302}.
\newblock


\bibitem[Usmanova et~al\mbox{.}(2024)]%
        {usmanova2024log}
\bibfield{author}{\bibinfo{person}{Ilnura Usmanova}, \bibinfo{person}{Yarden As}, \bibinfo{person}{Maryam Kamgarpour}, {and} \bibinfo{person}{Andreas Krause}.} \bibinfo{year}{2024}\natexlab{}.
\newblock \showarticletitle{Log barriers for safe black-box optimization with application to safe reinforcement learning}.
\newblock \bibinfo{journal}{\emph{Journal of Machine Learning Research}} \bibinfo{volume}{25}, \bibinfo{number}{171} (\bibinfo{year}{2024}), \bibinfo{pages}{1--54}.
\newblock


\bibitem[Wang et~al\mbox{.}(2020)]%
        {wang2020non}
\bibfield{author}{\bibinfo{person}{Allen Wang}, \bibinfo{person}{Ashkan Jasour}, {and} \bibinfo{person}{Brian~C Williams}.} \bibinfo{year}{2020}\natexlab{}.
\newblock \showarticletitle{Non-gaussian chance-constrained trajectory planning for autonomous vehicles under agent uncertainty}.
\newblock \bibinfo{journal}{\emph{IEEE Robotics and Automation Letters}} \bibinfo{volume}{5}, \bibinfo{number}{4} (\bibinfo{year}{2020}), \bibinfo{pages}{6041--6048}.
\newblock


\bibitem[Wang and Gros(2022)]%
        {wang2022solving}
\bibfield{author}{\bibinfo{person}{Kai Wang} {and} \bibinfo{person}{S{\'e}bastien Gros}.} \bibinfo{year}{2022}\natexlab{}.
\newblock \showarticletitle{Solving mission-wide chance-constrained optimal control using dynamic programming}. In \bibinfo{booktitle}{\emph{2022 IEEE 61st Conference on Decision and Control (CDC)}}. IEEE, \bibinfo{pages}{2947--2952}.
\newblock


\bibitem[Wang et~al\mbox{.}({[n.\,d.]})]%
        {wang2019neural}
\bibfield{author}{\bibinfo{person}{Lingxiao Wang}, \bibinfo{person}{Qi Cai}, \bibinfo{person}{Zhuoran Yang}, {and} \bibinfo{person}{Zhaoran Wang}.} \bibinfo{year}{[n.\,d.]}\natexlab{}.
\newblock \showarticletitle{Neural Policy Gradient Methods: Global Optimality and Rates of Convergence}. In \bibinfo{booktitle}{\emph{International Conference on Learning Representations}}.
\newblock


\bibitem[Wisniewski and Bujorianu(2023)]%
        {wisniewski2023probabilistic}
\bibfield{author}{\bibinfo{person}{Rafal Wisniewski} {and} \bibinfo{person}{Manuela~L Bujorianu}.} \bibinfo{year}{2023}\natexlab{}.
\newblock \showarticletitle{Probabilistic safety guarantees for Markov decision processes}.
\newblock \bibinfo{journal}{\emph{IEEE Trans. Automat. Control}} (\bibinfo{year}{2023}).
\newblock


\bibitem[Yuan et~al\mbox{.}(2023)]%
        {yuan2023linear}
\bibfield{author}{\bibinfo{person}{Rui Yuan}, \bibinfo{person}{Simon~Shaolei Du}, \bibinfo{person}{Robert~M. Gower}, \bibinfo{person}{Alessandro Lazaric}, {and} \bibinfo{person}{Lin Xiao}.} \bibinfo{year}{2023}\natexlab{}.
\newblock \showarticletitle{Linear Convergence of Natural Policy Gradient Methods with Log-Linear Policies}. In \bibinfo{booktitle}{\emph{The Eleventh International Conference on Learning Representations}}.
\newblock
\urldef\tempurl%
\url{https://openreview.net/forum?id=-z9hdsyUwVQ}
\showURL{%
\tempurl}


\bibitem[Yuan et~al\mbox{.}(2022)]%
        {yuan2022general}
\bibfield{author}{\bibinfo{person}{Rui Yuan}, \bibinfo{person}{Robert~M Gower}, {and} \bibinfo{person}{Alessandro Lazaric}.} \bibinfo{year}{2022}\natexlab{}.
\newblock \showarticletitle{A general sample complexity analysis of vanilla policy gradient}. In \bibinfo{booktitle}{\emph{International Conference on Artificial Intelligence and Statistics}}. PMLR, \bibinfo{pages}{3332--3380}.
\newblock


\end{thebibliography}
\appendix
\section{Proofs}
\subsection{Proof of Proposition \ref{pro_smooth}}\label{sec_pro_smooth}
\begin{proof}[Proof of Proposition \ref{pro_smooth}]
 For simplicity, we define the reward for a trajectory $\tau$ as 
\[r(\tau):=\sum_{t=0}^{H-1} r_t(s_t, a_t) + r_H(s_H).\]

To establish the Lipschitz continuity of $\rew$, we bound the norm of the gradient $\nabla_\theta \rew$. By policy gradient theorem~\cite[Theorem A.5]{kohler2017sub}, we have
 \begin{align}
     \nabla_\theta \rew
     = &\E^{\theta}\left[\left(\sum_{t\in[H-1]}\nabla_\theta\log \tpi^{\theta_t}_t(a_t|\ts_t)\right)r(\tau) \mid \ts_0 \right].
 \end{align}
 Since $\nabla_\theta \log\tpi^{\theta_t}_t(a|\ts)=(0,\dots,0,\nabla_{\theta_t} \log\tpi^{\theta_t}_t(a|\ts)^T,0,\dots, 0)^T$ for any $t\in[H-1]$, we have
 \begin{align}
     \|\nabla_\theta \rew\|
     \le &H\E^{\theta}\Bigl[\bigl\|\sum_{t\in[H-1]}\nabla_\theta\log \tpi^{\theta_t}_t(a_t|\ts_t)\bigr\|\mid \ts_0 \Bigr]\\
     = & H\E^{\theta}\Bigl[\sqrt{\sum_{t\in[H-1]}\left\|\nabla_{\theta_t}\log \tpi^{\theta_t}_t(a_t|\ts_t)\right\|^2} \mid \ts_0\Bigr]\\
     \le & M_gH^{\frac{3}{2}},\label{eq_boundvr}
 \end{align}
where we use Assumption~\ref{assumption_smooth_policy} in the last step.
To prove the smoothness of $\rew$, we bound the norm of the Hessian $\nabla_\theta^2 \rew$, which is calculated as follows:
\begin{small}
\begin{align}
&\nabla_\theta^2 \rew\\
= &\nabla_\theta\E^{\theta}\left[\left(\sum_{t\in[H-1]}\nabla_\theta\log \tpi^{\theta_t}_t(a_t|\ts_t)\right) r(\tau) \mid \ts_0\right] \\
= &\underbrace{\E^{\theta}\left[\left(\sum_{t\in[H-1]}\nabla_\theta\log \tpi^{\theta_t}_t(a_t|\ts_t)\right)\left(\sum_{t\in[H-1]}\nabla_\theta\log \tpi^{\theta_t}_t(a_t|\ts_t)\right)^T r(\tau) \mid \ts_0\right]}_{(i)}\\
&+\underbrace{\E^{\theta}\left[\left(\sum_{t\in[H-1]}\nabla_\theta^2\log \tpi^{\theta_t}_t(a_t|\ts_t)\right) r(\tau) \mid \ts_0\right]}_{(ii)}.
\end{align}
\end{small}
For term $(i)$, we have
\begin{align}
    \|(i)\|&\le H\, \E^{\theta}\left[\left\|\sum_{t\in[H-1]}\nabla_\theta\log \tpi^{\theta_t}_t(a_t|\ts_t)\right\|^2\mid \ts_0\right]\\
    &=H\sum_{t\in[H-1]}\E^{\theta}\left[\left\|\nabla_{\theta_t}\log \tpi^{\theta_t}_t(a_t|\ts_t)\right\|^2\mid \ts_0\right]\\
    &\le M_g^2 H^2 ,
\end{align}
where the last inequality uses Assumption \ref{assumption_smooth_policy}.
For term $(ii)$, we have
\begin{align}
    \|(ii)\|\le& H\E^{\theta}\left[\|\sum_{t\in[H-1]}\nabla_\theta^2\log \tpi^{\theta_t}_t(a_t|\ts_t)\|\mid \ts_0\right]\\
    =& H\max_{t\in[H-1]}\E^{\theta}\left[\|\nabla_\theta^2\log \tpi^{\theta_t}_t(a_t|\ts_t)\|\mid \ts_0\right]\\
    \le & M_hH.
\end{align}
where the last step follows from Assumption \ref{assumption_smooth_policy}. Combining the upper bounds for terms $(i)$ and $(ii)$, we bound \(\|\nabla_\theta^2 \rew\|\) as 
\begin{align}
    \|\nabla_\theta^2 \rew\|\le M_hH+M_g^2H^2.
\end{align}
To prove the Lipschitz continuity and smoothness constants for \(\cost\), we apply a similar analysis as above but replace  $r(\tau)$ with $c(\ts_H)$, which is upper bounded by 1. Consequently, the Lipschitz continuity and smoothness constants for $\cost$ are a factor of $H$ smaller than those for $\rew$.

For Point 2, since we use independent and identically distributed samples for $n$ trajectories, \(\ecost\) is unbiased by definition, and \(\hat\nabla_\theta \cost\) and \(\hat\nabla_\theta \rew\) are unbiased by the policy gradient theorem \cite[Theorem A.5]{kohler2017sub}.

To prove Point 3, we establish a concentration bound for the estimator $\ecost$ by Hoeffding’s inequality. For any $\varepsilon>0$, we have
\begin{align}
    \PP\left(\left|\ecost-\cost\right|\ge \varepsilon\right)\le 2\exp\left(-2n\varepsilon
    ^2\right). \nonumber
\end{align}
Rewriting this inequality, we obtain, for any $\beta\in(0,1)$,
\begin{align}
    \PP\left(\left|\ecost-\cost\right|\le \sigma_c^0(n)\sqrt{\ln{\frac{2}{\beta}}}\right)\ge 1-\beta,
    \end{align}
where $\sigma_c^0(n):=\frac{1}{\sqrt{2n}}$.

To prove Point 4, we first bound the norms of the estimators \(\hat\nabla_\theta \cost\) and \(\hat\nabla_\theta \rew\). For \(\hat\nabla_\theta \rew\), we have
\begin{align}
    &\left\|\hat\nabla_\theta \rew\right\|\\
    =&\Bigl\|\frac{1}{n}\!\!\sum_{i\in[n-1]}\!\sum_{t\in[H-1]}\!\!\!\nabla_\theta\log \tpi^{\theta_t}_t(a^i_t|\ts^i_t) (\!\sum_{m=t}^{H-1} {r}_m({s}^i_m, a^i_m) + {r}_H({s}^i_H) )\Bigr\|\\
    =&\sqrt{\sum_{t\in[H-1]}\Bigl\|\frac{1}{n}\sum_{i\in[n-1]}\nabla_{\theta_t}\log \tpi^{\theta_t}_t(a^i_t|\ts^i_t) (H-t)\Bigr\|^2}\\
    \le & M_g \sqrt{\int_1^{H} t^2dt}\le {M_gH^{\frac{3}{2}}}.
\end{align}
For \(\hat\nabla_\theta \cost\), we have 
\begin{align}
    &\bigl\|\hat\nabla_\theta \cost\bigr\|\\
    =&\Bigl\|\frac{1}{n}\sum_{i\in[n-1]}\sum_{t\in[H-1]}\nabla_\theta\log \tpi^{\theta_t}_t(a^i_t|\ts^i_t)){c}(\tilde{s}^i_H)\Bigr\|\\
    \le& \sqrt{\frac{1}{n}\sum_{i\in[n-1]}\sum_{t\in[H-1]}\Bigl\|\nabla_{\theta_t}\log \tpi^{\theta_t}_t(a^i_t|\ts^i_t))\Bigr\|^2}\\
    \le& H^{\frac{1}{2}}M_g.
\end{align}
Therefore, we can bound the error of the estimator \(\hat\nabla_\theta \rew\) as follows:
\begin{align}
    \|\hat\nabla_\theta \rew -\nabla_\theta \rew\|
    \le&\|\nabla_\theta \rew\|+ \|\nabla_\theta \rew\|\\
    \le& 2M_gH^{\frac{3}{2}},
\end{align}
where we use inequality \eqref{eq_boundvr} in the last step. The variance of the estimator \(\hat\nabla_\theta \rew\) can be upper bounded as follows:
\begin{align}
    &\E \|\hat\nabla_\theta \rew -\nabla_\theta \rew\|^2\\
    \le&\E\|\hat\nabla_\theta \rew\|^2 + \|\nabla_\theta \rew\|^2 + 2\|\nabla_\theta \rew\|\E\|\hat\nabla_\theta \rew\|\\
    \le& 4M_g^2H^3,
\end{align}
where we use inequality \eqref{eq_boundvr} in the last step. Similarly, for the estimator \(\hat\nabla_\theta \cost\), we have:
\begin{align}
    &\|\hat\nabla_\theta \cost -\nabla_\theta \cost\|\le 2 H^{\frac{1}{2}}M_g,\\
    &\E \|\hat\nabla_\theta \cost -\nabla_\theta \cost\|^2\le 4H M_g^2.
\end{align}
Therefore, we establish concentration bounds for the estimators {$\hat \nabla_\theta \rew$} and {$\hat \nabla_\theta \cost$} using the vector Bernstein inequality \cite[Lemma 18]{kohler2017sub}.  For any $\beta\in(0,1)$, if $n\ge 8\ln{\frac{\exp{\frac{1}{4}}}{\beta}}$, we have
\begin{align}  
&\PP\left(\left\|\hat \nabla_\theta \rew- \nabla_\theta \rew\right\|\le  \sigma_{r}^1(n)\sqrt{\ln{\frac{\exp{\frac{1}{4}}}{\beta}}}\right)
\ge 1-\beta,\nonumber
\end{align}
where $\sigma_{r}^1(n):=\frac{2\sqrt{2}L_r}{\sqrt{n}}$. Applying the same analysis as above, with $L_{r}$ replaced by $L_c$, we can establish the concentration bound for the estimator {$\hat \nabla_\theta \cost$}.
\end{proof}
\subsection{Proof of Theorem \ref{main}}\label{sec_pro_thm}
First, we prove that the iterates generated by the LB-SGD algorithm remain within a distance of $\Omega(\eta)$ from the boundary.
\begin{lemma}\label{small}
Let Assumptions~\ref{Assumption_continuity},~\ref{assumption_smooth_policy}, and \ref{emf} hold. Fix a confidence level $\beta\in(0,1)$. By setting {$\eta \leq \min\{\nu_{s},p\}$, $n = \mathcal{O}\left(\eta^{-2} \ln{\beta}^{-1}\right)$}, we have
\begin{align*}
    \PP\left\{\min_{i\in [I]} V_{c,0}^{\theta(i)}(\ts_0)-\delta\ge c\eta\right\}\ge 1-I\beta.
\end{align*}
with $c:=\frac{\ell}{24L_c}$.
\end{lemma}
\begin{proof}
The above results follow from \cite[Lemma 6]{usmanova2024log}.
\end{proof}

Next, we establish the gradient dominance property of the log barrier function, which allows us to bound the optimality of the policy {$\tpi_\theta$} given {$\nabla_\theta B_\eta^\theta$}.

\begin{lemma}\label{gdm}
Let Assumptions \ref{assumption_smooth_policy} and \ref{fn} hold. For any {\(\theta \in \Theta\),} we have
\begin{align*}
   V_{r,0}^{\tpi_M^*}(\ts_0)-\rew\le&\frac{H^{\frac{1}{2}}M_g}{\mu_F}\left\| \nabla_\theta B_\eta^{\theta}\right\|+ \!\eta\!+ \!\sqrt{\varepsilon_{bias}}H\left(\!1\!+\!\frac{\eta}{\cost-\delta}\right). 
\end{align*}
\end{lemma}

\begin{proof}
By performance difference lemma~\cite[Lemma A.3]{klein2024beyond}, we have
\begin{align}
  &V_{r,0}^{\tpi^*}(\ts_0)-\rew+\eta\frac{V_{c,0}^{\tpi^*}(\ts_0)-\cost}{ \cost}\\
  =& \sum_{t\in[H-1]} \mathbb{E}_{(\ts,a)\sim d_t^{\tpi^*}}  \left[A_{r,t}^{\theta}(\ts,a)+\frac{\eta
  A_{c,t}^{\theta}(\ts,a)}{\cost}\right].\label{s1}
\end{align}
Applying Jensen's inequality transfer error, we obtain the following for all $h\in[H-1]$ and $k \in\{r,c\}$:
\begin{align}
    \mathbb{E}_{(\ts,a)\sim d_t^{\tpi^*}}  \left[A_{k,t}^\theta(\ts,a)-{\mu_{k,t}}^T\nabla_{\theta_t}\log\tpi_t^{\theta_t}(a|\ts)\right]\le\sqrt{\varepsilon_{bias}}.\label{s2}
\end{align}
Plugging inequality \eqref{s2} into \eqref{s1}, we get
\begin{align}
    &V_{r,0}^{\tpi^*}(\ts_0)-\rew+\eta\frac{V_{c,0}^{\tpi}(\ts_0)-\cost}{ \cost}\\
  \le& \sum_{t\in[H-1]} \underbrace{\mathbb{E}_{(\ts,a)\sim d_t^{\tpi^*}} \left[\left({\mu_{r,t}}+\frac{\eta \mu_{c,t}}{ \cost}\right)^T \nabla_{\theta_t}\log\tpi_t^{\theta_t}(a|\ts)\right]}_{(t)}\\
  &+\sqrt{\varepsilon_{bias}}H\left(1+\frac{\eta}{\cost-\delta}\right).\label{ss1}
\end{align}
Next, we substitute the mathematical expressions for $\mu_{r,t}$ and $\mu_{k,t}$ into each term $(t)$:
\begin{align}
      (t)=&\E_{(\ts,a)\sim d_t^\theta}\Bigl[\nabla_{\theta_t} \log\tpi_t^{\theta_t}(a|\ts)^T (A_{r,t}^\theta(\ts,a)+\frac{A_{c,t}^\theta(\ts,a)}{\cost})\Bigr]^T \\
      &\cdot \left(F_t^\theta\right)^\dagger \mathbb{E}_{(\ts,a)\sim d_t^{\tpi^*}} \Bigl[\nabla_{\theta_t}\log\tpi_t^{\theta}(a|\ts)\Bigr].\label{eq_t}
\end{align}    
By the policy gradient theorem \cite[Theorem A.5]{kohler2017sub}, we have
\begin{align}
    &\nabla_\theta B_\eta^\theta \\
    =& \sum_{t\in[H-1]}\E_{(\ts,a)\sim d_t^\theta}\Bigl[\nabla_\theta \log\tpi_t^{\theta_t}(a|\ts)^T (A_{r,t}^\theta(\ts,a)+\frac{A_{c,t}^\theta(\ts,a)}{\cost})\Bigr]\\
    =&\begin{pmatrix}
\E_{(\ts,a)\sim d_0^\theta}\Bigl[\nabla_{\theta_0} \log\tpi_0^{\theta_0}(a|\ts)^T (A_{r,0}^\theta(\ts,a)+\frac{A_{c,0}^\theta(\ts,a)}{\cost})\Bigr]\\
\dots\\
\E_{(\ts,a)\sim d_{H-1}^\theta}\Bigl[\nabla_{\theta_{H-1}} \log\tpi_{H-1}^{\theta_{H-1}}(a|\ts)^T (A_{r,{H-1}}^\theta(\ts,a)+\frac{A_{c,{H-1}}^\theta(\ts,a)}{\cost})\Bigr]
\end{pmatrix}\label{eq_gradientb}
\end{align}
We define the matrices $F$ and $L$ as
\begin{equation*}
F:= \begin{bmatrix}
F_0^\theta & \hdots & 0\\
\vdots   &  \ddots    & \vdots \\
0 & \hdots & F_{H-1}^\theta
\end{bmatrix}, L:= \begin{pmatrix}
\mathbb{E}_{(\ts,a)\sim d_0^{\tpi^*}} \Bigl[\nabla_{\theta_0}\log\tpi_0^{\theta_0}(a|\ts)\Bigr]\\
\dots\\
\mathbb{E}_{(\ts,a)\sim d_{H-1}^{\tpi^*}} \Bigl[\nabla_{\theta_{H-1}}\log\tpi_{H-1}^{\theta_{H-1}}(a|\ts)\Bigr]
\end{pmatrix}.
\end{equation*}
Using the above notation along with \eqref{eq_gradientb} on \eqref{eq_t}, we find
\begin{align}
    \sum_{t\in[H-1]}(t)&= \nabla_\theta B_\eta^\theta F^\dagger L\\
    &=\left( \mathbf{P}_{\textbf{Im}(F)}\nabla_\theta B_\eta^\theta\right) F^\dagger  L \\
    &\le \frac{H^{\frac{1}{2}}M_g}{\mu_F}\|\mathbf{P}_{\textbf{Im}(F)}\nabla_\theta B_\eta^\theta\|\\
    &\le \frac{H^{\frac{1}{2}}M_g}{\mu_F}\|\nabla_\theta B_\eta^\theta\|,
\end{align}
where $\mathbf{P}_{\textbf{Im}(F)}\nabla_\theta B_\eta^\theta$ denotes the orthogonal projection onto the image space $F$, and we use Assumptions \ref{assumption_smooth_policy} and \ref{fn} in the first inequality. Substituting the above inequality into \eqref{ss1}, we obtain
\begin{align*}
   V_{r,0}^{\tpi^*}(\ts_0)-\rew\le&\frac{H^{\frac{1}{2}}M_g}{\mu_F}\left\| \nabla_\theta B_\eta^{\theta}\right\|+ \!\eta\!+ \!\sqrt{\varepsilon_{bias}}H\left(\!1\!+\!\frac{\eta}{\cost-\delta}\right). 
\end{align*}
\end{proof}
Using the above results, we proceed with the proof of Theorem \ref{main}.
\begin{proof}[Proof of Theorem \ref{main}]
By setting the sample size $n$ as
\begin{align}
    n=\mathcal{O}\left(\eta^{-4} \allowbreak\ln{\beta}^{-1}\right),
\end{align}
Point (1) of Theorem \ref{main} is satisfied with a probability of at least \( 1 - I\beta \), as established by \cite[Lemma 6]{usmanova2024log}. With this choice of \( n \), we achieve a lower bound on the distance of the iterates from the boundary, and the following holds almost surely:

\begin{enumerate}
    \item By \cite[Lemma 1]{usmanova2024log}, \(\|\Delta(i)\| \le \frac{\eta}{4}\), where \(\Delta(i) := \hat\nabla_\theta B_\eta^{\theta(i)} - \nabla_\theta B_\eta^{\theta(i)}\).
    \item By \cite[Lemma 5]{usmanova2024log}, \(\min_{i\in[I]} \gamma(i) \ge C\eta\), where \(C\) is defined in Theorem \ref{main}.
    \item By \cite[Lemma 2]{usmanova2024log}, \(\gamma(i) \leq \frac{1}{M(i)}\), where \(M(i)\) is the local smoothness constant of the log barrier function \(B_\eta^{\theta(i)}\).
\end{enumerate}
With the above information established, we proceed to prove the second point of Theorem \ref{main}. Since \(\gamma(i) \leq \frac{1}{M(i)}\), we can upper bound \(B_\eta^{\theta(i+1)} - B_\eta^{\theta(i)}\) as follows:
\begin{align}
    &B_\eta^{\theta(i+1)} - B_\eta^{\theta(i)} \\
    \geq & \gamma(i) \left\langle \nabla_\theta B_\eta^{\theta(i)}, \hat \nabla_\theta B_\eta^{\theta(i)} \right\rangle - \frac{M(i) \gamma(i)^2}{2} \left\|\hat \nabla_\theta B_\eta^{\theta(i)}\right\|^2 \\
    \geq & \gamma(i) \left\langle \nabla_\theta B_\eta^{\theta(i)}, \hat \nabla_\theta B_\eta^{\theta(i)} \right\rangle - \frac{\gamma(i)}{2} \left\|\hat \nabla_\theta B_\eta^{\theta(i)}\right\|^2 \\
    = & \gamma(i) \left\langle \nabla_\theta B_\eta^{\theta(i)}, \Delta(i) + \nabla_\theta B_\eta^{\theta(i)} \right\rangle - \frac{\gamma(i)}{2} \left\|\Delta(i) + \nabla_\theta B_\eta^{\theta(i)}\right\|^2 \\
    = & \frac{\gamma(i)}{2} \left\| \nabla_\theta B_\eta^{\theta(i)} \right\|^2 - \frac{\gamma(i)}{2} \left\|\Delta(i)\right\|^2.\label{ss_1}
\end{align}

We divide the analysis into two cases based on the \textbf{if condition} in Algorithm \ref{alg:cap}, line 4.

\textbf{Case 1:} If \(\|\hat \nabla_\theta B_\eta^{\theta(i)}\| \ge \frac{\eta}{2}\), then \(\|\nabla_\theta B_\eta^{\theta(i)}\| \ge \frac{\eta}{4}\) since \(\|\Delta(i)\| \le \frac{\eta}{4}\). Thus, we rewrite \eqref{ss_1} as:
\begin{align}
    B_\eta^{\theta(i+1)} - B_\eta^{\theta(i)} &\ge \frac{C\eta^2}{8} \left\|\nabla_\theta B_\eta^{\theta(i)}\right\| - \frac{C\eta^3}{32},\label{case1:1}
\end{align}
where we used \(\gamma(i) \ge C\eta\) in the last step. Since \(\min_{i\in[I]} V_{c,0}^{\theta(i)}(\ts_0) \geq c\eta\), we can express Lemma \ref{gdm} as:
\begin{align}
    V_{r,0}^{\tpi^*}(\ts_0) - V_{r,0}^{\theta(i)}(\ts_0) 
    \le a + \frac{ H^{\frac{1}{2}}M_g}{\mu_F} \left\| \nabla_\theta B_\eta^{\theta(i)} \right\|,\label{gmdp}
\end{align}
where \(a := \eta + {\sqrt{\varepsilon_{bias}}}H\left(1 + \frac{1}{c}\right)\). Plugging \eqref{gmdp} into \eqref{case1:1}, we have
\begin{align}
    B_\eta^{\theta(i+1)} - B_\eta^{\theta(i)} 
    \ge & \frac{C\mu_F\eta^2}{8 H^{\frac{1}{2}}M_g} \left(V_{r,0}^{\tpi^*}(\ts_0) - V_{r,0}^{\theta(i)}(\ts_0)\right) - \frac{C\eta^3}{32} - \frac{aC\mu_F\eta^2 }{8 H^{\frac{1}{2}}M_g}.\nonumber
\end{align}
This can be further simplified to:
\begin{align}
    V_{r,0}^{\tpi^*}(\ts_0) - V_{r,0}^{\theta(i+1)}(\ts_0) 
    \le & \left(1 - \frac{C\mu_F\eta^2}{8 H^{\frac{1}{2}}}M_g\right)\left(V_{r,0}^{\tpi^*}(\ts_0) - V_{r,0}^{\theta(i)}(\ts_0)\right) + \frac{C\eta^3}{32} \nonumber \\
    & + \frac{aC\mu_F\eta^2 }{8 H^{\frac{1}{2}}M_g} + \eta \log\frac{V_{c,0}^{\theta(i)}(\ts_0)-\delta}{V_{c,0}^{\theta(i)}(\ts_0)-\delta}.\nonumber
\end{align}

By recursively applying the above inequality and ensuring that \(\frac{C\mu_F\eta^2}{8 H^{\frac{1}{2}}M_g} < 1\), we obtain:
\begin{align}
        &V_{r,0}^{\tpi^*}(\ts_0) - V_{r,0}^{\theta(i+1)}(\ts_0)\\
        \le&\left(1- \frac{C\mu_F\eta^2}{8H^{\frac{1}{2}}M_g}\right)^{i+1}\left(V_{r,0}^{\tpi^*}(\ts_0) - V_{r,0}^{\theta(0)}(\ts_0)\right)+\underbrace{\eta\log (V_{c,0}^{\theta(i+1)}(\ts_0)-\delta)}_{(i)}\\
        &-\underbrace{\eta\left(1- \frac{C\mu_F\eta^2}{8H^{\frac{1}{2}}M_g}\right)^{i}\log (V_{c,0}^{\theta(0)}(\ts_0)-\delta)}_{(ii)}\\
        &+\underbrace{\left(\frac{C\eta^3}{32}  + \frac{aC\mu_F\eta^2 }{8H^{\frac{1}{2}}M_g}\right)\sum_{j\in[i]}\left(1- \frac{C\mu_F\eta^2}{8H^{\frac{1}{2}}M_g}\right)^j}_{(iii)}\\
        & -\underbrace{\frac{C\mu_F\eta^3}{8H^{\frac{1}{2}}M_g}\sum_{j=1}^{i}\left(1- \frac{C\mu_F\eta^2}{8H^{\frac{1}{2}}M_g}\right)^{i-j} \log(V_{c,0}^{\theta(j)(\ts_0)}-\delta)}_{(iv)}.\label{eq_all}
    \end{align}
We analyze each term in \eqref{eq_all}. For term (i): Since \(\cost \leq 1\), we can upper bound this term:
   \begin{align}
       (i) \leq \eta \log(1 - \delta).
   \end{align}
For term (ii): Given that \(V_{c,0}^{\theta(0)(\ts_0)} - \delta \geq \nu_s\) by Point (3) of Assumption \ref{Assumption_continuity}, we can lower bound this term:
   \begin{align}
       (ii) \geq \eta \left(1 - \frac{C\mu_F\eta^2}{8H^{\frac{1}{2}}M_g}\right)^{i} \log \nu_s \geq \eta \log \nu_s.
   \end{align}
For term (iii): We can upper bound it as:
   \begin{align}
       (iii) \leq \left(\frac{C\eta^3}{32} + \frac{aC\mu_F\eta^2}{8H^{\frac{1}{2}}M_g}\right) \left(1 - \left(1 - \frac{C\mu_F\eta^2}{8H^{\frac{1}{2}}M_g}\right)\right)^{-1} = \frac{ H^{\frac{1}{2}}\eta M_g}{4\mu_F} + a.
   \end{align}
For term (iv): Since \(V_{c,0}^{\theta(j)}(\ts_0) - \delta \geq c\eta\) for all \(j \in [i]\) by Lemma \ref{small}, we can lower bound this term as:
   \begin{align}
       (iv) \geq \frac{C\mu_F\eta^3}{8H^{\frac{1}{2}}M_g} \sum_{j=1}^{i} \left(1 - \frac{C\mu_F\eta^2}{8H^{\frac{1}{2}}M_g}\right)^{i-j} \log(c\eta) \geq \eta \log(c\eta).
   \end{align}
Combining the above, we can write \eqref{eq_all} as:
\begin{align}
    V_{r,0}^{\tpi^*}(\ts_0) - V_{r,0}^{\theta(i+1)}(\ts_0) 
    \leq& \left(1 - \frac{C\mu_F\eta^2}{8H^{\frac{1}{2}}M_g}\right)^{i+1} \left(V_{r,0}^{\tpi^*}(\ts_0) - V_{r,0}^{\theta(0)}(\ts_0)\right) \nonumber \\
    & + \frac{ H^{\frac{1}{2}}\eta M_g}{4\mu_F} + a + \eta \log \frac{1 - \delta}{c\eta\nu_s}. \label{m1}
\end{align}

\textbf{Case 2:} If \(\|\hat \nabla_\theta B_\eta^{\theta(i)}\| \leq \frac{\eta}{2}\), we have \(\|\nabla_\theta B_\eta^{\theta(i)}\| \leq \|\hat \nabla_\theta B_\eta^{\theta(i)}\| + \|\Delta(i)\| \leq \frac{3\eta}{4}\). By \eqref{gmdp}, we have:
\begin{align}
    V_{r,0}^{\tpi^*}(\ts_0) - V_{r,0}^{\theta(i)}(\ts_0) 
    \leq a + \frac{ 3H^{\frac{1}{2}}\eta M_g}{4\mu_F}. \label{m2}
\end{align}

Combining the inequalities \eqref{m1} and \eqref{m2}, we conclude the proof for Point 2 of Theorem \ref{main}.
\end{proof}
\end{document}